\documentclass[times,review,sort&compress]{elsarticle}
\usepackage{epstopdf}
\usepackage[caption=false]{subfig}
\usepackage[hidelinks]{hyperref}
\usepackage[active]{srcltx}
\usepackage{amsthm}
\usepackage{amssymb}
\usepackage{stackrel}
\usepackage{amsmath,amscd}

\newtheorem{dl}{Theorem}[section]
\newtheorem{tl}[dl]{Corollary}
\newtheorem{yl}[dl]{Lemma}
\newtheorem{dy}[dl]{Definition}
\newtheorem{lz}[dl]{Example}

\newtheorem{remark}[dl]{Remark}

\numberwithin{equation}{section}

\newcommand{\be}{\begin{equation}}
\newcommand{\ee}{\end{equation}}
\newcommand{\ba}{\begin{array}}
\newcommand{\ea}{\end{array}}
\newcommand{\bmn}{\begin{eqnarray}}
\newcommand{\emn}{\end{eqnarray}}
\newcommand{\bnm}{\begin{eqnarray*}}
\newcommand{\enm}{\end{eqnarray*}}
\newcommand{\bln}{\begin{subequations}}
\newcommand{\eln}{\end{subequations}}

\newcommand{\poq}[2]{(#1;q)_{#2}}

\def\qed{\hfill \rule{4pt}{7pt}}
\def\pf{\noindent {\it Proof.} }

\begin{document}

\title{General $q$-series transformations based on Abel's lemma on  summation by parts and their applications}
\author{Jianan Xu\fnref{fn1}}
\fntext[fn1]{E-mail address: 20224007007@stu.suda.edu.cn}
\address[P.R.China]{Department of Mathematics, Soochow University, SuZhou 215006, P.R.China}
\author{Xinrong Ma\fnref{fn3}}
\fntext[fn3]{Corresponding author. E-mail address: xrma@suda.edu.cn.}
\address[P.R.China]{Department of Mathematics, Soochow University, SuZhou 215006, P.R.China}
\markboth{J. Xu and  X. Ma}{General $q$-series transformations based on Abel's lemma on summation by parts and their applications}
\begin{abstract}
In this   paper,  we  establish three new and general transformations with sixteen parameters and bases via Abel's lemma on summation by parts. As applications, we set up a lot of new transformations of basic hypergeometric series.  Among include some new results  of Gasper and Rahman's quadratic,
cubic, and quartic transformations. Furthermore, we put forward  the so-called $(R,S)$-type transformation with arbitrary degree  to unify such  multibasic transformations. Some special $(R,S)$-type transformations are presented.\end{abstract}
\maketitle
\vspace{20pt}
\parskip 7pt
\section{Introduction}
It is well known  that Abel's lemma on summation by parts is one of the most  important tools in  Analysis and Number Theory, which is analogous to integration by parts. Abel's lemma can often be stated as follows:
\begin{yl}[Abel's lemma on summation by parts: {\rm \cite[p. 313]{knopp}}]\label{abelemma} For integer $n\geq 0$ and sequences
$\{A_n\}_{n\geq 0}$ and $\{B_n\}_{n\geq 0}$, it always holds
\begin{align}\sum_{k=0}^{n-1} B_k\Delta A_k
=A_{n-1}B_{n}-A_{-1}B_0+\sum_{k=0}^{n-1}A_k\nabla B_k.\label{abelparts}
\end{align}
Hereafter, we define the backward and forward difference operators $\Delta$ and $\nabla$ acting on arbitrary sequence $\{X_n\}_{n\geq 0}$, respectively, by
\begin{align}
\Delta X_k:=X_k-X_{k-1}\quad and \quad
\nabla X_k:=X_k-X_{k+1}.
\end{align}
\end{yl}
Simple as it seems, Abel's lemma has been proved to be one of basic and useful tools in  various infinite or finite ($q$-) series. The reader is referred to \cite{bha,cam,chon,kom,jht,zhang}  for this topic.  It is of interest to note that in \cite{chenthree, chenthree-1,bian}, Chen and his collaborators associated Abel's  lemma with the well-known WZ-algorithm. Especially worthy of mention is that since the 1990s,  Chu and his collaborators have systematically investigated and demonstrated many transformation and summation formulae from the theory of numerous (general, basic, and theta) hypergeometric functions using Abel's lemma on summation by parts. See \cite{chulatest,chu06,chu06c,chu2007,chu2008,Chu1,Wang}. A careful look at  Chu et al.'s  elementary proofs  reveals that there are two basic algebraic identities  underlying many of their arguments in loc. cit.
\begin{yl}\label{lemm12}For any complex parameters $a,b,c,x,y,z,w: xyzw\neq 0$, the modified Jacobi theta function is defined to be
\begin{align}&\theta(x;q):=(x;q)_\infty(q/x;q)_\infty,~~\mathcal{D}(x):=1-x.\label{charafunchi}
\end{align}Then the following two identities are true.
\begin{enumerate} [(i)]
\item (Four-term algebraic identity)
 \begin{align}
 \mathcal{D}\left(cx,\frac{x}{c},bz,\frac{z}{b}\right)-\mathcal{D}\left(bx,\frac{x}{b},cz,\frac{z}{c}\right)=\frac{z}{c}\mathcal{D}\left(bc,\frac{c}{b},xz,\frac{x}{z}\right).\label{kkklll-111}
 \end{align}
\item (Weierstrass' theta identity: \cite[Ex. 20.5.6]{10})
 \begin{align}
 \theta\left(cx,\frac{x}{c},bz,\frac{z}{b};q\right)-\theta\left(bx,\frac{x}{b},cz,\frac{z}{c};q\right)=\frac{z}{c}\theta\left(bc,\frac{c}{b},xz,\frac{x}{z};q\right).\label{trivalweierstrass}
 \end{align}
\end{enumerate}
In the above, the multivariate notation
$\mathcal{D}(x_1,x_2,\ldots,x_m):\:=\:\mathcal{D}(x_1)\mathcal{D}(x_2)\ldots\mathcal{D}(x_m),$ $m\geq 1$ being integer and
\(\theta(x_1,x_2,\ldots,x_m):\:=\:\theta(x_1)\theta(x_2)\ldots\theta(x_m).\)
\end{yl}
It is worth mentioning that   \eqref{kkklll-111}  is equivalent to Weierstrass' theta identity \eqref{trivalweierstrass}. As for this equivalency, we refer the reader to \cite{wangjin} for  a full proof. In addition, we would like to refer the reader to \cite{koornwinder} for a short history on the origins of Weierstrass' theta identity. As far as its applications to theta and elliptic hypergeometric series are concerned,  the reader might consult \cite{liu,liu-1} by Liu, \cite{xrma} by Ma, \cite{rosen} by Rosengren and Schlosser, and \cite{warnaartheta} by Warnaar for further details.

Apart from Chu et al.'s works, it is worthy of attention that in their papers \cite{Gasper89, Gasper90}, Gasper and Rahman et al. applied Abel's lemma on summation by parts, although  not stated clearly,  to set up  a lot of multibasic (quadratic, cubic, and quartic) transformations for $q$-series. More results afterward can be found in \cite{rahman} by Rahman and Verma, as well as \cite{jainverma} by  Jain and Verma.  What intrigued us is  that in \cite{Gasper89}, Gasper believed that their method could not work in a general setting, the reason is that for an arbitrary sequence
\begin{align}
s_n=\frac{(a ; p)_n(b ; q)_n(c ; P)_n(d ; Q)_n}{\left(e ; p^{\prime}\right)_n\left(f ; q^{\prime}\right)_n\left(g; P^{\prime}\right)_n\left(h ; Q^{\prime}\right)_n},
\end{align}
 the corresponding difference $\Delta(s_n)=s_n-s_{n-1}$ does not take the factorial  form. The reader might consult \cite[Eq. (2.8)]{Gasper89}
for further  details. Unlike his intuition,  there do exist some  general sequences $s_n$ such that $\Delta(s_n)$ can be decomposed into product of factorial  factors.   This finding prompted us
to make  a systematic study of Gasper's multibasic transformations. To our knowledge, Gasper's  results were extended to the setting of root systems like the $U(n)$-form  \cite{bha-1} given by Bhatanagar and Miner, \cite{rosen} by Rosengren and Schlosser. Other generalizations can be found in \cite{roselin,zhangys}.

In the present paper, we will   focus on  applications of \eqref{kkklll-111} and Abel's lemma on summation by parts to the theory of basic hypergeometric series.

Our paper is organized as follows. In the next   section, we will  show three general eight-base transformations, i.e., Theorems \ref{type-i}, \ref{type-ii}, and \ref{type-iii}  by Abel's lemma and \eqref{kkklll-111}.  Section 3 is devoted to some concrete  transformations that are deducible from  Theorems \ref{type-i} and \ref{type-ii}. Among include  some new  $q$-series transformations that seem  peculiar.  In Section 4, based on a general but more concrete form of Theorem \ref{type-i},  we will  propose the so-called $(R,S)$-type transformation with arbitrary degree, in order to better understand Gasper-Rahman's research works on  multibasic transformations.  Some special $(R,S)$-type transformations are presented.

{\bf Notation:} Throughout this paper, we will  follow the notation and terminology in the book \cite{10} by Gasper and Rahman. The $q$-shifted factorial of a complex variable $a$ with the base $q,|q|<1$, is given by
\begin{align*}
(a;q)_{\infty}:=\prod_{k=0}^{\infty}(1-aq^k), ~~(a;q)_{n}:=\frac{(a;q)_{\infty}}{(aq^n;q)_{\infty}}~(n\in \mathbf{Z}),
\end{align*}
where $\mathbf{Z}$ denotes the set of integers.
For  any complex numbers $a_1,a_2,\ldots,a_m$, the multivariate notation
\begin{align*}
(a_1,a_2,\ldots,a_m;q)_n:=\prod_{k=1}^{m}(a_k;q)_n
\end{align*}
and the basic hypergeometric series with the base $q$ and the variable $z$ is defined to be
\begin{align}
_{r}\phi_{r-1}\left[\begin{matrix}a_1,a_2,\ldots,a_r\\ b_1,b_2,\ldots,b_{r-1}\end{matrix};q,z\right]
:=\sum_{n=0}^{+\infty}\frac{(a_1,a_2,\ldots,a_r;q)_n}{(b_1,b_2,\ldots,b_{r-1};q)_n}\frac{z^{n}}{\poq{q}{n}}.\label{phiseries-1}
\end{align}
\section{The main theorems and proofs}
\subsection{The first transformation}
By virtue  of the basic algebraic identity \eqref{kkklll-111}, it is easy to establish
\begin{dl}[The first transformation]\label{type-i}
For any  complex vectors $\bar{a}=(a_1,a_2,a_3,a_4),$
$\bar{p}=(p_1,p_2,p_3,p_4)$, $\bar{x}=(x_1,x_2,x_3,x_4), \bar{q}=(q_1,q_2,q_3,q_4)$ with  $a_ix_ip_iq_i\neq 0,1\leq i\leq 4$, there holds
\begin{align}
&\sum_{k=0}^{n-1}\Gamma_{k-1}[\bar{a};\bar{p}]\frac{KL^{k-1}-1}{KL^{k-1} }\times \prod_{i=1}^4\frac{(a_i^2;p_i)_{k-1}}{(K/a_i^2;L/p_i)_k
}\prod_{i=1}^4
\frac{(x_i^2;q_i)_k}
{(K_0/x_i^2;L_0/q_i)_k}\nonumber\\
  &=
\prod_{i=1}^4\frac{(a_i^2;p_i)_{n-1}}
{(K/a_i^2;L/p_i)_{n-1}}\prod_{i=1}^4\frac{(x_i^2;q_i)_n}
{(K_0/x_i^2;L_0/q_i)_{n}}-
\prod_{i=1}^4\frac
{L/p_i-K/a_i^2}{p_i-a_i^2}\label{pppppp}\\
&+\sum_{k=0}^{n-1}\Gamma_{k}[\bar{x};\bar{q}]\frac{1-K_0L_0^{k}}{K_0L_0^{k} }
\times \prod_{i=1}^4\frac{(x_i^2;q_i)_k}
{(K_0/x_i^2;L_0/q_i)_{k+1}}
\prod_{i=1}^4\frac{(a_i^2;p_i)_k}
{(K/a_i^2;L/p_i)_k}\nonumber
.
\end{align}
Hereafter, we define
\begin{align}K:=a_1a_2a_3a_4,\quad& L:=\sqrt{p_1p_2p_3p_4};\label{kkklll-1}\\
K_0:=x_1x_2x_3x_4,\quad&L_0:=\sqrt{q_1q_2q_3q_4};\label{kkklll-2}\end{align}
and  the function
\begin{align}
\Gamma_k[\bar{x};\bar{q}]&:= \left(x_1x_2 (q_1q_2)^{\frac{k}{2}}-x_3x_4 (q_3q_4)^{\frac{k}{2}}\right) \label{chachacha}\\
&\times\left(x_1x_3 (q_1q_3)^{\frac{k}{2}}-x_2x_4
(q_2q_4)^{\frac{k}{2}}\right)\left(x_1x_4 (q_1q_4)^{\frac{k}{2}}-x_2x_3 (q_2q_3)^{\frac{k}{2}}\right).\nonumber
\end{align}
\end{dl}
\pf First of all, we choose two sequence~$\{A_k,B_k\}_{k\geq 0}$~in Lemma \ref{abelparts} by
\begin{subequations}
\begin{align}
A_k&:=\prod_{i=1}^4\frac{(a_i^2;p_i)_k}
{(K/a_i^2;L/p_i)_k};\label{ak}\\
B_k&:=\prod_{i=1}^4\frac{(x_i^2;q_i)_k}
{(K_0/x_i^2;L_0/q_i)_k}.\label{bk}
\end{align}
\end{subequations}
Now we compute the difference
\bnm
\Delta A_k=\bigg(1-\frac{A_{k-1}}{A_k}\bigg)A_k=\bigg(1-\frac{1}{K^2 L^2}\prod_{i=1}^4\frac{\left(K L-a_i^2 p_i\right)}{\left(1-a_i^2 {p_i}\right) }\bigg)\bigg|_{p_i\to p_i^{k-1},L\to L^{k-1}} A_k.
\enm
Herein and in what follows, $\left.F\right|_\sigma$ means applying the parameter  substitution $\sigma$ to the function $F$.
Since,  by \eqref{kkklll-111},
\begin{align*}
1-\frac{1}{K^2 L^2}\prod_{i=1}^4\frac{\left(K L-a_i^2 p_i\right)}{\left(1-a_i^2 {p_i}\right) }&=\frac{KL-1}{KL}\frac{ \left(a_1a_2 \sqrt{p_1p_2}-a_3a_4\sqrt{p_3p_4}\right)}{\left(1-a_1^2p_1\right)\left(1-a_2^2 p_2\right) }\\
&\times\frac{ \left(a_1a_3\sqrt{p_1p_3}-p_2p_4 \sqrt{p_2p_4}\right)\left(a_1a_4 \sqrt{p_1p_4}-a_2a_3 \sqrt{p_2p_3}\right) }{ \left(1-a_3^2 p_3\right)\left(1-a_4^2 p_4\right)},
\end{align*}
 reformulating the last identity  in terms of $\Gamma_{k-1}[\bar{a};\bar{p}]$, we have
\begin{align}
\Delta A_k=\Gamma_{k-1}[\bar{a};\bar{p}]\frac{KL^{k-1}-1}{KL^{k-1} }\prod_{i=1}^4\frac{(a_i^2;p_i)_{k-1}}{(K/a_i^2;L/p_i)_k}.\label{diffak}
\end{align}
In the meantime, we compute directly
\bnm
\nabla B_k=\bigg(1-\frac{B_{k+1}}{B_k}\bigg)B_k=\bigg(1-K_0^2 L_0^2 \prod_{i=1}^4\frac{\left(1-x_i^2 {q_i}\right)}{\left(K_0 L_0-x_i^2 q_i\right)}\bigg)\bigg|_{q_i\to q_i^{k},L_0\to L_0^{k}}\times B_k.
\enm
By \eqref{kkklll-111}, we easily find
\bnm
&&1-K_0^2 L_0^2 \prod_{i=1}^4\frac{\left(1-x_i^2 {q_i}\right)}{\left(K_0 L_0-x_i^2 q_i\right)}=
K_0L_0\left(1-K_0L_0\right)\\
&&\quad\times\frac{\left(x_1x_2\sqrt{q_1q_2}-x_3x_4 \sqrt{q_3q_4}\right) \left(x_1x_3\sqrt{q_1q_3}-x_2x_4 \sqrt{q_2q_4}\right) \left(x_1x_4 \sqrt{q_1q_4}-x_2x_3\sqrt{q_2q_3}\right)}{\left(K_0L_0-x_1^2q_1\right)\left(K_0L_0-x_2^2 q_2\right)  \left(K_0L_0-x_3^2 q_3\right)\left(K_0L_0-x_4^2 q_4\right)},
\enm
which enables us to reformulate $\nabla B_k$   as the form
\begin{align}
\nabla B_k=\Gamma_{k}[\bar{x};\bar{q}]\frac{1-K_0L_0^{k}}{K_0L_0^{k} }
\times \prod_{i=1}^4\frac{(x_i^2;q_i)_k}
{(K_0/x_i^2;L_0/q_i)_{k+1}}.\label{diffbk}
\end{align}
Consequently, we substitute the differences \eqref{diffak} and \eqref{diffbk} into \eqref{abelparts},  finding
\begin{align*}
&\sum_{k=0}^{n-1}\Gamma_{k-1}[\bar{a};\bar{p}]\frac{KL^{k-1}-1}{KL^{k-1} }\times \prod_{i=1}^4\frac{(a_i^2;p_i)_{k-1}}{(K/a_i^2;L/p_i)_k
}\prod_{i=1}^4
\frac{(x_i^2;q_i)_k}
{(K_0/x_i^2;L_0/q_i)_k}\\
  &=A_{n-1}B_n-A_{-1}B_0\\
&+\sum_{k=0}^{n-1}\Gamma_{k}[\bar{x};\bar{q}]\frac{1-K_0L_0^{k}}{K_0L_0^{k} }
\times \prod_{i=1}^4\frac{(x_i^2;q_i)_k}
{(K_0/x_i^2;L_0/q_i)_{k+1}}
\prod_{i=1}^4\frac{(a_i^2;p_i)_k}
{(K/a_i^2;L/p_i)_k}
.
\end{align*}
Recall (cf. \cite[(I.10)]{10})
 that\begin{align}
(a ; q)_{n-k}=\frac{(a ; q)_n}{\left(q^{1-n} / a ; q\right)_k}\left(-\frac{q}{a}\right)^k q^{k(k-1)/2-n k}.\label{vvvv}
\end{align}
Thus it is easy to compute
\[A_{-1}=
\prod_{i=1}^4\frac{(a_i^2;p_i)_{-1}}
{(K/a_i^2;L/p_i)_{-1}}=
\prod_{i=1}^4\frac
{L/p_i-K/a_i^2}{p_i-a_i^2}.\]
Finally, the desired transformation follows from a direct substitution of $A_n, A_{-1}$, $B_n$ and a bit of simplification.
\qed

\subsection{The second transformation}
In a very similar way, we can apply the following equivalent form of \eqref{kkklll-111}
\begin{align}
 \mathcal{D}\left(cx,\frac{x}{c},bz,\frac{z}{b}\right)-\frac{z}{c}\mathcal{D}\left(bc,\frac{c}{b},xz,\frac{x}{z}\right)=\mathcal{D}\left(bx,\frac{x}{b},cz,\frac{z}{c}\right)\label{kkklll-111-new}
 \end{align}
and Abel's lemma on summation by parts to show
\begin{dl}[The second transformation]\label{type-ii}With the same notation as Theorem \ref{type-i}. For any integer $n$, it holds
\begin{align}
&\sum_{k=0}^{n-1}(-1)^k\frac{\Gamma_{k-1}[\bar{a};\bar{p}]}
{(K_0;L_0)_k}\frac{KL^{k-1}-1}{KL^{k-1} }~ \prod_{i=1}^4\frac{(a_i^2;p_i)_{k-1}(x_i^2;q_i)_k}{(K/a_i^2;L/p_i)_k
}\prod_{j=0}^{k-1}
\frac{K_0L_0^{j}}{\Gamma_j[\bar{x};\bar{q}]}\nonumber\\
  &=\frac{(-1)^n}{(K_0;L_0)_n}\prod_{i=1}^4\frac{(a_i^2;p_i)_{n-1}(x_i^2;q_i)_n}
{(K/a_i^2;L/p_i)_{n-1}} \prod_{j=0}^{n-1}
\frac{K_0L_0^j}{\Gamma_j[\bar{x};\bar{q}]}-\prod_{i=1}^4
\frac{L/p_i-K/a_i^2}{p_i-a_i^2}\label{identityfinal-II}\\
&+\sum_{k=0}^{n-1}\frac{(-1)^k}
{(K_0;L_0)_{k+1}}\prod_{j=0}^{k}
\frac{K_0L_0^j}{\Gamma_j[\bar{x};\bar{q}]}
\prod_{i=1}^4\frac{(a_i^2;p_i)_k(x_i^2;q_i)_k\left(1-K_0L_0^k/(x_i^2q_i^k)\right)}
{(K/a_i^2;L/p_i)_k}.\nonumber
\end{align}
\end{dl}
\pf It suffices to choose $A_k$~and~$B_k$~in Abel’s lemma \ref{abelparts} respectively to be
\begin{align}
A_k&:=\prod_{i=1}^4\frac{(a_i^2;p_i)_k}
{(K/a_i^2;L/p_i)_k};\label{aak}\\
B_k&:=\frac{(-1)^k}{(K_0;L_0)_k} \prod_{j=0}^{k-1}
\frac{K_0L_0^j}{\Gamma_j[\bar{x};\bar{q}]}\prod_{i=1}^4(x_i^2;q_i)_k.\label{bbk}
\end{align}
In the preceding theorem,  we have already evaluated $
\Delta A_k$ to be given by \eqref{diffak}. Thus, all that remains is to find the difference
\bnm
\nabla B_k=\bigg(1-\frac{B_{k+1}}{B_k}\bigg)~B_k=\bigg(1-\frac{1}{(K_0L_0^k-1)}
\frac{K_0L_0^k}{\Gamma_k[\bar{x};\bar{q}]}\prod_{i=1}^4(1-x_i^2q_i^{k})\bigg)~B_k.
\enm
Note that from \eqref{kkklll-111-new}, it follows
\bnm
1-\frac{1}{(K_0L_0^k-1)}
\frac{K_0L_0^k}{\Gamma_k[\bar{x};\bar{q}]}\prod_{i=1}^4(1-x_i^2q_i^{k})=\frac{1}{K_0L_0^k\left(1-K_0L_0^k\right)}\frac{1}{\Gamma_k[\bar{x};\bar{q}]}\prod_{i=1}^4\left(K_0L_0^k-x_i^2q_i^k\right).
\enm
This gives rise to
\begin{align}
\nabla B_k=\frac{(-1)^k}{K_0^2L_0^{2k}(K_0;L_0)_{k+1}} \prod_{j=0}^{k}
\frac{K_0L_0^j}{\Gamma_j[\bar{x};\bar{q}]}\prod_{i=1}^4\left(K_0L_0^k-x_i^2q_i^k\right)(x_i^2;q_i)_k.\label{diffbbk}
\end{align}
By substituting \eqref{diffak} and \eqref{diffbbk} into \eqref{abelparts}, we obtain
\begin{align*}
&\sum_{k=0}^{n-1}(-1)^k\frac{\Gamma_{k-1}[\bar{a};\bar{p}]}
{(K_0;L_0)_k}\frac{KL^{k-1}-1}{KL^{k-1} }\prod_{j=0}^{k-1}
\frac{K_0L_0^{j}}{\Gamma_j[\bar{x};\bar{q}]}~ \prod_{i=1}^4\frac{(a_i^2;p_i)_{k-1}(x_i^2;q_i)_k}{(K/a_i^2;L/p_i)_k
}\\
  &=A_{n-1}B_n-A_{-1}B_0\\
&+\frac{1}
{K_0^2}\sum_{k=0}^{n-1}\frac{(-1)^k}
{L_0^{2k}(K_0;L_0)_{k+1}}\prod_{j=0}^{k}
\frac{K_0L_0^j}{\Gamma_j[\bar{x};\bar{q}]}
\prod_{i=1}^4\frac{(a_i^2;p_i)_k(x_i^2;q_i)_k\left(K_0L_0^k-x_i^2 q_i^k\right)}
{(K/a_i^2;L/p_i)_k}.
\end{align*} Finally, substitute $A_n$ and $B_n$ and make some simplifications by using the basic relation
\[\prod_{i=1}^4\left(K_0L_0^k-x_i^2 q_i^k\right)=K_0^2L_0^{2k}\prod_{i=1}^4\left(1-K_0L_0^k/x_i^2 q_i^k\right).\]  Hence \eqref{identityfinal-II} is proved.
\qed

\subsection{The third transformation}
To show this result, we need the following basic fact \cite[(2.9)]{Gasper89} used by Gasper.
\begin{yl}For any complex sequences $\{\alpha_n,\beta_n\}_{n\geq 0}$, there holds
\begin{align}
\sum_{k=0}^n\alpha_k\sum_{i=0}^{n-k}\beta_i=\sum_{k=0}^n\beta_k\sum_{i=0}^{n-k}\alpha_i.\label{alpha-beta}
\end{align}
\end{yl}
On applying \eqref{alpha-beta} to $A_n$ and $B_n$ of Theorem \ref{type-i}, we easily find that
\begin{dl}[The third transformation]\label{type-iii}With the same notation as Theorem \ref{type-i}. Then, for any integer $n\geq 0$, it holds
\begin{align}
\prod_{i=1}^4\frac{(x_i^2;q_i)_{n+1}}{(K_0/x_i^2;L_0/q_i)_{n+1}}&\bigg(-1+\sum_{k=1}^n\Gamma_{k-1}[\bar{a};\bar{p}]\frac{1-KL^{k-1}}{KL^{k-1} }\nonumber\\ &\qquad\times \prod_{i=1}^4\frac{(a_i^2;p_i)_{k-1}}{(K/a_i^2;L/p_i)_k}\frac{((L_0/q_i)^{-n}x_i^2/K_0;L_0/q_i)_{k}}{(q_i^{-n}/x_i^2;q_i)_{k}}\bigg)
\nonumber\\
=\prod_{i=1}^4\frac{(a_i^2;p_i)_{n}}{(K/a_i^2;L/p_i)_{n}}&\bigg(-1+\sum_{k=0}^n\Gamma_{k}[\bar{x};\bar{q}]\frac{1-K_0L_0^{k}}{K_0L_0^{k} }\label{alpha-beta-new}
\\
&\qquad\times \prod_{i=1}^4\frac{(x_i^2;q_i)_k}
{(K_0/x_i^2;L_0/q_i)_{k+1}} \frac{((L/p_i)^{1-n}a_i^2/K;L/p_i)_{k}}{(p_i^{1-n}/a_i^2;p_i)_{k}}\bigg).\nonumber
\end{align}
\end{dl}
\pf  As mentioned above,  we need only to specialize \eqref{alpha-beta} by
\begin{align*}
\alpha_k&=\Gamma_{k-1}[\bar{a};\bar{p}]\frac{KL^{k-1}-1}{KL^{k-1} }\prod_{i=1}^4\frac{(a_i^2;p_i)_{k-1}}{(K/a_i^2;L/p_i)_k
}=A_k-A_{k-1};\\
\beta_k&=\Gamma_{k}[\bar{x};\bar{q}]\frac{1-K_0L_0^{k}}{K_0L_0^{k} }
 \prod_{i=1}^4\frac{(x_i^2;q_i)_k}
{(K_0/x_i^2;L_0/q_i)_{k+1}}=B_k-B_{k+1},
\end{align*}
where $A_k$ and $B_k$ are given by \eqref{ak} and  \eqref{bk}, respectively.
It is clear that
\begin{align}\sum_{k=1}^n\alpha_kB_{n+1-k}
+\sum_{k=0}^n\beta_kA_{n-k}=A_n-B_{n+1},\label{alpha-beta-newnew}
\end{align}
which amounts to
\begin{align}
&\sum_{k=1}^n\Gamma_{k-1}[\bar{a};\bar{p}]\frac{KL^{k-1}-1}{KL^{k-1} }\prod_{i=1}^4\frac{(a_i^2;p_i)_{k-1}}{(K/a_i^2;L/p_i)_k}
\prod_{i=1}^4\frac{(x_i^2;q_i)_{n+1-k}}{(K_0/x_i^2;L_0/q_i)_{n+1-k}}\nonumber\\
&+\sum_{k=0}^n\Gamma_{k}[\bar{x};\bar{q}]\frac{1-K_0L_0^{k}}{K_0L_0^{k} }
 \prod_{i=1}^4\frac{(x_i^2;q_i)_k}
{(K_0/x_i^2;L_0/q_i)_{k+1}}
\prod_{i=1}^4\frac{(a_i^2;p_i)_{n-k}}
{(K/a_i^2;L/p_i)_{n-k}}=A_n-B_{n+1}.\label{alpha-beta-newnewnew}
\end{align}
Referring to \eqref{vvvv}, we easily compute\begin{align*}
&\prod_{i=1}^4\frac{(a_i^2;p_i)_{n-k}}{(K/a_i^2;L/p_i)_{n-k}}=\prod_{i=1}^4\frac{(a_i^2;p_i)_{n}}{(K/a_i^2;L/p_i)_{n}}\prod_{i=1}^4 \frac{((L/p_i)^{1-n}a_i^2/K;L/p_i)_{k}}{(p_i^{1-n}/a_i^2;p_i)_{k}};\\
&\prod_{i=1}^4\frac{(x_i^2;q_i)_{n+1-k}}{(K_0/x_i^2;L_0/q_i)_{n+1-k}}=\prod_{i=1}^4\frac{(x_i^2;q_i)_{n+1}}{(K_0/x_i^2;L_0/q_i)_{n+1}}\prod_{i=1}^4 \frac{((L_0/q_i)^{-n}x_i^2/K_0;L_0/q_i)_{k}}{(q_i^{-n}/x_i^2;q_i)_{k}}.\end{align*}
Upon substituting all these computations in \eqref{alpha-beta-newnewnew}, we  obtain  \eqref{alpha-beta-new} at once.
\qed

\section{Applications of the first and second transformations}This section is totally devoted to applications of Theorems \ref{type-i}, \ref{type-ii}, and \ref{type-iii}.
\subsection{Some concrete results deduced from the first transformation}
\subsubsection{Some $q$-series transformations covered by Theorem \ref{type-i}}
First of all, we proceed to show a more concrete result embraced by Theorem \ref{type-i}.
\begin{tl}With the same notation as Theorem \ref{type-i}. Then, for any integer $n\geq 0$, it holds
\begin{align}
&\sum_{k=0}^{n-1}\Gamma_{k}[\bar{x};\bar{q}]\bigg(\prod_{i=1}^4(1-K_0/x_i^2(L_0/q_i)^{k})+\prod_{i=1}^4(1-x_i^2q_i^k)\bigg)\frac{K_0L_0^{k}-1}{K_0L_0^{k}}\nonumber\\
&\qquad\times \bigg(\prod_{i=1}^4\frac{(x_i^2;q_i)_{k}}{(K_0/x_i^2;L/q_i)_{k+1}
}\bigg)^2=
\bigg(\prod_{i=1}^4\frac{(x_i^2;q_i)_{n}}
{(K_0/x_i^2;L_0/q_i)_{n}}\bigg)^2-1\label{pppppp-1-finalnew-2}.
\end{align}
\end{tl}
\pf  To establish \eqref{pppppp-1-finalnew-2}, it needs only to set   in \eqref{pppppp} that
\begin{align*}
(p_1,p_2,p_3,p_4)&=(q_1,q_2,q_3,q_4)\qquad\mbox{and}
\\
(a_1,a_2,a_3,a_4)&=(x_1q_1^{1/2},x_2q_2^{1/2},x_3q_3^{1/2},x_4q_4^{1/2}).
\end{align*}
As such, it is easily found that
\begin{align*}L=L_0,~~& K=K_0L_0;~~ \Gamma_{k-1}[\bar{a};\bar{p}]=\Gamma_{k}[\bar{x};\bar{q}];\\
(K/a_i^2;L/p_i)_k&=(K_0L_0/x_i^2q_i;L_0/q_i)_k=\frac{(K_0/x_i^2;L_0/q_i)_{k+1}}{1-K_0/x_i^2}.
\end{align*}
Hence, all these reduces  \eqref{pppppp} to
\begin{align}
&\sum_{k=0}^{n-1}\Gamma_{k}[\bar{x};\bar{q}]\frac{K_0L_0^{k}-1}{K_0L_0^{k}}\prod_{i=1}^4\frac{(x_i^2q_i;q_i)_{k-1}}{(K_0L_0/x_i^2q_i;L_0/q_i)_k
}\prod_{i=1}^4
\frac{(x_i^2;q_i)_k}
{(K_0/x_i^2;L_0/q_i)_k}\nonumber\\
  &=
\prod_{i=1}^4\frac{(x_i^2q_i;q_i)_{n-1}}
{(K_0L_0/x_i^2q_i;L_0/q_i)_{n-1}}\prod_{i=1}^4\frac{(x_i^2;q_i)_n}
{(K_0/x_i^2;L_0/q_i)_{n}}-
\prod_{i=1}^4\frac
{1-K_0/x_i^2}{1-x_i^2}\label{pppppp-1-final}\\
&+\sum_{k=0}^{n-1}\Gamma_{k}[\bar{x};\bar{q}]\frac{1-K_0L_0^{k}}{K_0L_0^{k} }
 \prod_{i=1}^4\frac{(x_i^2;q_i)_k}
{(K_0/x_i^2;L_0/q_i)_{k+1}}
\prod_{i=1}^4\frac{(x_i^2q_i;q_i)_k}
{(K_0L_0/x_i^2q_i;L_0/q_i)_k}\nonumber
.
\end{align}
By multiplying both sides of \eqref{pppppp-1-final} with
$\prod_{i=1}^4(1-x_i^2)/(1-K_0/x_i^2),$ we have
\begin{align}
&\sum_{k=0}^{n-1}\Gamma_{k}[\bar{x};\bar{q}]\frac{K_0L_0^{k}-1}{K_0L_0^{k}}\prod_{i=1}^4\frac{(x_i^2;q_i)_{k}}{(K_0/x_i^2;L_0/q_i)_{k+1}
}\prod_{i=1}^4
\frac{(x_i^2;q_i)_k}
{(K_0/x_i^2;L_0/q_i)_k}\nonumber\\
  &=
\prod_{i=1}^4\frac{(x_i^2;q_i)_{n}}
{(K_0/x_i^2;L_0/q_i)_{n}}\prod_{i=1}^4\frac{(x_i^2;q_i)_n}
{(K_0/x_i^2;L_0/q_i)_{n}}-1\nonumber\\
&+\sum_{k=0}^{n-1}\Gamma_{k}[\bar{x};\bar{q}]\frac{1-K_0L_0^{k}}{K_0L_0^{k} }
 \prod_{i=1}^4\frac{(x_i^2;q_i)_k}
{(K_0/x_i^2;L_0/q_i)_{k+1}}
\prod_{i=1}^4\frac{(x_i^2;q_i)_{k+1}}
{(K_0/x_i^2;L_0/q_i)_{k+1}}\nonumber
.
\end{align}
The conclusion is proved.
\qed

We now turn to one of the most important  cases of Theorem \ref{type-i}.
\begin{tl}For any integer $n\geq 0$, it holds
\begin{align}
&a_0\sum_{k=0}^{n-1}\frac{(1-bcdep^{2k}/a)p^{k}}{1-bcde/a}
\frac{(b,d,e,bc^2de/a^2;p)_{k}}
{(ap/c,bcep/a,bcdp/a,cdep/a;p)_{k}
}
\frac{(yq,wq,uq,yz^2wuq/x^2;q)_k}
{(xq/z,zwuq/x,yzuq/x,yzwq/x;q)_k}\nonumber\\
  &\qquad\qquad=\frac{(b,d,e,bc^2de/a^2;p)_{n}}
{(a/c,bce/a,bcd/a,cde/a;p)_{n}}\frac{(yq,wq,uq,yz^2wuq/x^2;q)_n}
{(xq/z,zwuq/x,yzuq/x,yzwq/x;q)_{n}}-1\label{rogerspsi65}
\\
&-a_1\sum_{k=1}^{n}\frac{(1-yzwuq^{2k}/x)q^{k}}{1-yzwu/x}\frac{(y,w,u,yz^2wu/x^2;q)_k}
{(xq/z,zwuq/x,yzuq/x,yzwq/x;q)_{k}}
\frac{(b,d,e,bc^2de/a^2;p)_{k}}
{(a/c,bce/a,bcd/a,cde/a;p)_{k}}\nonumber
,
\end{align}
where two coefficients
\begin{align}
a_0&=\frac{a(1-bc/a)(1-cd/a)(1-ce/a)(1-bcde/a)}{c(1-bcd/a)(1-bce/a)(1-cde/a)(1-a/c)};\label{chi0}\\
a_1&=\frac{x(1-yz/x)(1-zw/x)(1-zu/x)(1-yzwu/x)}{z(1-y)(1-w)(1-u)(1-yz^2wu/x^2)}.\label{chi1}
\end{align}
\end{tl}
\pf It suffices in \eqref{pppppp} to put
\begin{align*}
p_i=p, ~~(1\leq i\leq 4),~&
(a_1^2,a_2^2,a_3^2,a_4^2)=(bp,dp,ep,bc^2dep/a^2);\\
 q_i=q,~~(1\leq i\leq 4),~&(x_1^2,x_2^2,x_3^2,x_4^2)=(yq,wq,uq,yz^2wuq/x^2).
\end{align*}
These yield
\begin{align*}&K:=bcdep^2/a,~~L=p^2,\\
&\Gamma_{k-1}[\bar{a};\bar{p}]=-bde(1-bc/a)(1-cd/a)(1-ce/a)p^{3k};\\
&K_0:=yzwuq^2/x,~~L_0=q^2,\\
&\Gamma_{k}[\bar{x};\bar{q}]=-ywu(1-yz/x)(1-zw/x)(1-zu/x)q^{3k+3}.
\end{align*}
So we are able to reduce \eqref{pppppp} to
\begin{align*}
&-\frac{a(1-bc/a)(1-cd/a)(1-ce/a)}{c}\sum_{k=0}^{n-1}(bcdep^{2k}/a-1)p^{k}\\
&\times
\frac{(bp,dp,ep,bc^2dep/a^2;p)_{k-1}}
{(ap/c,bcep/a,bcdp/a,cdep/a;p)_k
}
\frac{(yq,wq,uq,yz^2wuq/x^2;q)_k}
{(xq/z,zwuq/x,yzuq/x,yzwq/x;q)_k}\nonumber\\
  &=\frac{(bp,dp,ep,bc^2dep/a^2;p)_{n-1}}
{(ap/c,bcep/a,bcdp/a,cdep/a;p)_{n-1}}\frac{(yq,wq,uq,yz^2wuq/x^2;q)_n}
{(xq/z,zwuq/x,yzuq/x,yzwq/x;q)_{n}}\\
&-\frac{(1-bcd/a)(1-bce/a)(1-cde/a)(1-a/c)}{(1-b)(1-d)(1-e)(1-bc^2de/a^2)}\\
&-\frac{xq(1-yz/x)(1-zw/x)(1-zu/x)}{z}\sum_{k=0}^{n-1}(1-yzwuq^{2k+2}/x)q^{k}\\
&\times \frac{(yq,wq,uq,yz^2wuq/x^2;q)_k}
{(xq/z,zwuq/x,yzuq/x,yzwq/x;q)_{k+1}}
\frac{(bp,dp,ep,bc^2dep/a^2;p)_k}
{(ap/c,bcep/a,bcdp/a,cdep/a;p)_k}.\nonumber
\end{align*}
As the last step, we multiply both sides of this identity with
\[\frac{(1-b)(1-d)(1-e)(1-bc^2de/a^2)}{(1-bcd/a)(1-bce/a)(1-cde/a)(1-a/c)},\]
 obtaining\begin{align*}
&-\frac{a(1-bc/a)(1-cd/a)(1-ce/a)}{c(1-bcd/a)(1-bce/a)(1-cde/a)(1-a/c)}\sum_{k=0}^{n-1}(bcdep^{2k}/a-1)p^{k}\\
&\times
\frac{(b,d,e,bc^2de/a^2;p)_{k}}
{(ap/c,bcep/a,bcdp/a,cdep/a;p)_{k}
}
\frac{(yq,wq,uq,yz^2wuq/x^2;q)_k}
{(xq/z,zwuq/x,yzuq/x,yzwq/x;q)_k}\nonumber\\
  &=\frac{(b,d,e,bc^2de/a^2;p)_{n}}
{(a/c,bce/a,bcd/a,cde/a;p)_{n}}\frac{(yq,wq,uq,yz^2wuq/x^2;q)_n}
{(xq/z,zwuq/x,yzuq/x,yzwq/x;q)_{n}}-1
\\
&-\frac{x(1-yz/x)(1-zw/x)(1-zu/x)}{z(1-y)(1-w)(1-u)(1-yz^2wu/x^2)}\sum_{k=1}^{n}(1-yzwuq^{2k}/x)q^{k}\\
&\times \frac{(y,w,u,yz^2wu/x^2;q)_k}
{(xq/z,zwuq/x,yzuq/x,yzwq/x;q)_{k}}
\frac{(b,d,e,bc^2de/a^2;p)_{k}}
{(a/c,bce/a,bcd/a,cde/a;p)_{k}}.\nonumber
\end{align*}
 Rewritten in terms of $a_0,a_1$, then desired result comes out.
\qed\\

Transformation \eqref{rogerspsi65} contains the following  results of interest as special cases, some of which we believe are new in the literature,  like  \eqref{important-3.4}--\eqref{1.39}.
\begin{lz}[{\rm \cite[Ex. 2.5]{10}}] For any integer $n\geq 0$, we have
\begin{align}
 \sum_{k=0}^{n}\frac{1-ywuq^{2k}}{1-ywu}\frac{(y,w,u,ywu;q)_k}
{(q,wuq,yuq,ywq;q)_{k}}q^{k} = \frac{(yq,wq,uq,ywuq;q)_n}
{(q,wuq,yuq,ywq;q)_{n}}.\end{align}
\end{lz}
\pf  It only needs to set $a = bc$ and $x = z$ in \eqref{rogerspsi65}.
\qed

Moreover, letting  $w=x$ and then taking $x\to 0$ in  \eqref{rogerspsi65}, we have immediately
\begin{lz}Let $a_0$ be given by \eqref{chi0}. Then, for any integer $n\geq 0$, we have
\begin{align}
&a_0\sum_{k=0}^{n-1}\frac{1-bcdep^{2k}/a}{1-bcde/a}
\frac{(b,d,e,bc^2de/a^2;p)_{k}}
{(ap/c,bcep/a,bcdp/a,cdep/a;p)_{k}
}
\frac{(yq,uq;q)_k}
{(zuq,yzq;q)_k}(pz)^k\nonumber\\
  &=\frac{(b,d,e,bc^2de/a^2;p)_{n}}
{(a/c,bce/a,bcd/a,cde/a;p)_{n}}\frac{(yq,uq;q)_n}
{(zuq,yzq;q)_{n}}z^n-\frac{(1-y z) (1-u z)}{z(1-y) (1-u) }\label{important-3.4}
\\
&+\frac{(1-z)(1-yzu)}{z(1-y)(1-u)}\sum_{k=0}^{n}\frac{1-yzuq^{2k}}{1-yzu} \frac{(y,u;q)_k}
{(zuq,yzq;q)_{k}}
\frac{(b,d,e,bc^2de/a^2;p)_{k}}
{(a/c,bce/a,bcd/a,cde/a;p)_{k}}z^k\nonumber
.\end{align}
\end{lz}

Furthermore, the limiting case $n\to\infty$ of  \eqref{important-3.4} yields
\begin{lz}Assume  $|z|<1$. Then it holds
\begin{align}
&a_0\sum_{k=0}^{\infty}\frac{1-bcdep^{2k}/a}{1-bcde/a}
\frac{(b,d,e,bc^2de/a^2;p)_{k}}
{(ap/c,bcep/a,bcdp/a,cdep/a;p)_{k}
}
\frac{(yq,uq;q)_k}
{(zuq,yzq;q)_k}(pz)^k\nonumber\\
  &=-\frac{(u z-1) (y z-1)}{z(u-1) (y-1) }\label{magic}
\\
&+\frac{(1-z)(1-yzu)}{z(1-y)(1-u)}\sum_{k=0}^{\infty}\frac{1-yzuq^{2k}}{1-yzu} \frac{(y,u;q)_k}
{(zuq,yzq;q)_{k}}
\frac{(b,d,e,bc^2de/a^2;p)_{k}}
{(a/c,bce/a,bcd/a,cde/a;p)_{k}}z^k\nonumber
.\end{align}
\end{lz}

In the case of  $p=q$, dividing both sides of \eqref{magic} by $a_0$ and  taking the limit of  $c\to a$, we obtain
\begin{lz}Assume that $|z|<1$. We have
\begin{align}
&\sum_{k=0}^{\infty}\frac{1-bde q^{2k}}{1-bde}
\frac{(b,d,e,bde,yq,uq;q)_{k}}
{(q,beq,bdq,deq,zuq,yzq;q)_{k}
}(qz)^k\\
  &=\frac{(1-z)(1-yzuq^2)}{(1-zuq)(1-yzq)}\sum_{k=0}^{\infty}\frac{1-yzuq^{2k+2}}{1-yzuq^2}
\frac{(yq,uq,bq,dq,eq,bdeq;q)_{k}}
{(q, zuq^2,yzq^2,beq,bdq,deq;q)_{k}}z^{k}.\nonumber
\end{align}
Or equivalently,
\begin{align}&{}_8\phi_{7}\left[\begin{matrix}bde,&q(bde)^{1/2},-q(bde)^{1/2},b,d,e,yq,uq\\ &(bde)^{1/2},-(bde)^{1/2},beq,bdq,deq,zuq,yzq\end{matrix};q,qz\right]\label{phiseries-2}\\
&=\frac{(1-z)(1-yzuq^2)}{(1-zuq)(1-yzq)}~
{}_8\phi_{7}\left[\begin{matrix}bdeq,&q^2(yzu)^{1/2},-q^2(yzu)^{1/2},bq,dq,eq,yq,uq\\ &q(yzu)^{1/2},-q(yzu)^{1/2},beq,bdq,deq,zuq^2,yzq^2\end{matrix};q,z\right].\nonumber
\end{align}
\end{lz}

Alternatingly, letting  $e\to 0$ in \eqref{rogerspsi65},  we find at once
 \begin{lz}
\begin{align}
&\frac{a(1-bc/a)(1-cd/a)}{c(1-bcd/a)(1-a/c)}\sum_{k=0}^{n-1}
\frac{(b,d;p)_{k}}
{(ap/c,bcdp/a;p)_{k}
}p^{k}~
\frac{(yq,wq,uq,yz^2wuq/x^2;q)_k}
{(xq/z,zwuq/x,yzuq/x,yzwq/x;q)_k}\nonumber\\
  &=\frac{(b,d;p)_{n}}
{(a/c,bcd/a;p)_{n}}\frac{(yq,wq,uq,yz^2wuq/x^2;q)_n}
{(xq/z,zwuq/x,yzuq/x,yzwq/x;q)_{n}}-1\label{important-3.7}
\\
&-a_1\sum_{k=1}^{n}\frac{(1-yzwuq^{2k}/x)q^{k}}{1-yzwu/x} \frac{(y,w,u,yz^2wu/x^2;q)_k}
{(xq/z,zwuq/x,yzuq/x,yzwq/x;q)_{k}}
\frac{(b,d;p)_{k}}
{(a/c,bcd/a;p)_{k}}\nonumber
.
\end{align}
\end{lz}
Moreover, on taking  $b\to \infty$ in the last identity and applying the algebraic identity \eqref{kkklll-111} and \eqref{vvvv}, we come up with
\begin{lz}
\begin{align}
&\frac{a(1-cd/a)}{cd(1-a/c)}\sum_{k=0}^{n-1}
\frac{(d;p)_{k}}
{(ap/c;p)_{k}
}\bigg(\frac{a}{cd}\bigg)^k
\frac{(yq,wq,uq,yz^2wuq/x^2;q)_k}
{(xq/z,zwuq/x,yzuq/x,yzwq/x;q)_k}\nonumber\\
  &=\frac{(d;p)_{n}}
{(a/c;p)_{n}}\bigg(\frac{a}{cd}\bigg)^n\frac{(yq,wq,uq,yz^2wuq/x^2;q)_n}
{(xq/z,zwuq/x,yzuq/x,yzwq/x;q)_{n}}\\
&\qquad\quad-\frac{(1-x/z)(1-ywz/x)(1-yzu/x)(1-wzu/x)}{(1-y)(1-w)(1-u)(1-yz^2wu/x^2)}\nonumber
\\
&-a_1\sum_{k=0}^{n}\frac{1-yzwuq^{2k}/x}{1-yzwu/x}\frac{(y,w,u,yz^2wu/x^2;q)_k}
{(xq/z,zwuq/x,yzuq/x,yzwq/x;q)_{k}}
\frac{(d;p)_{k}}
{(a/c;p)_{k}}\bigg(\frac{aq}{cd}\bigg)^k\nonumber
.
\end{align}
\end{lz}
Under the assumption that $|a/cd|<1$,  it is reasonable to  take $n\to \infty$. So we have
\begin{lz}For $|a/cd|<1$, it holds
\begin{align}
&\frac{a(1-cd/a)}{cd(1-a/c)}\sum_{k=0}^{\infty}
\frac{(d;p)_{k}}
{(ap/c;p)_{k}
}\bigg(\frac{a}{cd}\bigg)^k
\frac{(yq,wq,uq,yz^2wuq/x^2;q)_k}
{(xq/z,zwuq/x,yzuq/x,yzwq/x;q)_k}\nonumber\\
  &=-\frac{(1-x/z)(1-ywz/x)(1-yzu/x)(1-wzu/x)}{(1-y)(1-w)(1-u)(1-yz^2wu/x^2)}\label{1.39}
\\
&-a_1\sum_{k=0}^{\infty}\frac{1-yzwuq^{2k}/x}{1-yzwu/x}\frac{(y,w,u,yz^2wu/x^2;q)_k}
{(xq/z,zwuq/x,yzuq/x,yzwq/x;q)_{k}}
\frac{(d;p)_{k}}
{(a/c;p)_{k}}\bigg(\frac{aq}{cd}\bigg)^k.\nonumber
\end{align}
\end{lz}

Now we consider the special case of $x=z$ and $p=q$. In this case, we find that \eqref{1.39} can be expressed in the
standard notation of basic hypergeometric series as follows.
\begin{lz} For $|a/cd|<1$, it holds
\begin{align}&{}_5\phi_{4}\left[\begin{matrix}yq,wq,uq,ywuq,d\\ wuq,yuq,ywq,aq/c\end{matrix};q,\frac{a}{cd}\right]\label{phiseries-1-new}\\
&=\frac{1-a/c}{1-a/cd}~
{}_7\phi_{6}\left[\begin{matrix}ywu,&q(ywu)^{1/2},-q(ywu)^{1/2},y,w,u,d\\ &(ywu)^{1/2},-(ywu)^{1/2},yuq,wuq,ywq,a/c\end{matrix};q,\frac{aq}{cd}\right].\nonumber
\end{align}
\end{lz}
\subsubsection{Some new results on Gasper's quadratic, cubic, and quartic transformations}
To continue with applications of  Theorem \ref{type-i},  we now consider the quadratic transformation studied by Gasper \cite[Ex. 3.26]{10}. By applying Theorem \ref{type-i}, we find a new recursive relation.
\begin{dl} \label{thm3.2}Let us define
\begin{align}
G(a):=\sum_{k=0}^{\infty}\left(1-a b c q^{3k-1}\right)
\frac{(a b c/q,d, q/ d; q)_{k}\left(a c/q, ab/q,bcq; q^2\right)_{k}}
{\left(q^2,a b c q/ d, a b c d; q^2\right)_{k}(a/q,b q, c q ; q)_{k}}
 q^k.\label{ooooo}
\end{align}
Then it holds
\begin{align}
G(a)&=\frac{(abc/q;q)_{\infty}(a/d;q^2)_{\infty}(ad/q;q^2)_{\infty}}{(a/q;q)_{\infty}(abcq/d;q^2)_{\infty}(abcd;q^2)_{\infty}}G(0)\label{gasperid-new}\\
&-a~
\frac{(abc/q,d, q/ d ; q)_{\infty}}{\left(q^2,a b c q/ d, a b c d; q^2\right)_{\infty}} \frac{\left(a c q, abq,bcq; q^2\right)_{\infty}}{(a/q,b q, c q ; q)_{\infty}}\sum_{k=0}^\infty\frac{(a/d,ad/q;q^2)_{k}}{\left(ac q, abq; q^2\right)_{k}}q^{2k-1}.\nonumber
\end{align}
\end{dl}
\pf It suffices to take  in \eqref{pppppp} that
\begin{align*}
p_1=p_2=q,~~&p_3=p_4=q^2,~~(a_1^2,a_2^2,a_3^2,a_4^2)=(dq,q^2/d, bcq^3,a^2bcq^2);\\
 q_1=q_2=q,~~&q_3=q_4=q^2,~~(x_1^2,x_2^2,x_3^2,x_4^2)=(abc,q^2/a,acq, abq).
\end{align*}
It is easy to check
\begin{align*}K:=abcq^4,&~~L=q^3,~~\Gamma_{k-1}[\bar{a};\bar{p}]=bc(1-a/d)(ad-q)(1-abcq^k)q^{4k+1};\\
K_0:=abcq^2,~~&L_0=q^3,~\Gamma_{k}[\bar{x};\bar{q}]=bc(q-ab)(q-ac)(1-aq^k)q^{4k+2}.
\end{align*}
Thus \eqref{pppppp} can be simplified as
\begin{align}
&\frac{(1-a/d)(q-ad)}{a}\sum_{k=0}^{n-1}(1-abcq^k)\left(1-a b c q^{3 k+1}\right) q^k\nonumber\\
&\qquad\qquad\times\frac{(dq, q^2/ d ; q)_{k-1}
\left(b c q^3, a^2bcq^2 ; q^2\right)_{k-1}}{\left(a b c q^3 / d, a b c d q^2 ; q^2\right)_k(a q,q^2/a ; q)_k}
 \frac{(a b c,q^2/a; q)_k\left(a c q, abq; q^2\right)_k}
{\left(q^2,  a^2 b c; q^2\right)_k(b q, c q ; q)_k} \label{xinrong-today}
\\
&=\frac{(dq, q^2/ d ; q)_{n-1}
\left(b c q^3, a^2bcq^2 ; q^2\right)_{n-1}}{\left(a b c q^3 / d, a b c d q^2 ; q^2\right)_{n-1}(a q,q^2/a ; q)_{n-1}}\times \frac{(a b c,q^2/a; q)_{n}
\left(a c q, abq; q^2\right)_{n}}{\left(q^2,  a^2 b c; q^2\right)_{n}(b q, c q ; q)_{n}}-A_{-1} \nonumber\\
&+\frac{(q-ab)(q-ac)}{a}\sum_{k=0}^{n-1}(1-aq^k)\left(1-a b c q^{3 k+2}\right)q^k
 \nonumber\\
&\qquad\qquad\times\frac{(a b c; q)_k\left(a c q, abq; q^2\right)_k}
{\left(q^2,  a^2 b c; q^2\right)_{k+1}(b q, c q ; q)_{k+1}} \frac{(dq, q^2/ d ; q)_{k}
\left(b c q^3, a^2bcq^2 ; q^2\right)_{k}}{\left(a b c q^3 / d, a b c d q^2 ; q^2\right)_k(a q; q)_k}.\nonumber
\end{align}
Upon multiplying both sides of \eqref{xinrong-today} by
\[\frac{a(1-d)(1-q/d)(1-bcq)(1-a^2bc)}{(1-abc)(1-a/d)(q-ad)},\] we obtain
\begin{align*}
&\sum_{k=0}^{n-1}\left(1-a b c q^{3k+1}\right) \frac{(abcq, d, q/ d ; q)_{k}
\left(b c q,acq,abq; q^2\right)_{k}}{\left(q^2, a b c q^3 / d, a b c d q^2 ; q^2\right)_k(aq,b q, c q ; q)_k} q^k \\
&=\frac{a(1-q^{n+1}/a)(1-aq^n)}{(1-q^{2n})(1-abc)(1-a/d)(q-ad)}\frac{(abc,d, q/ d ; q)_{n}
}{\left(q^2,a b c q^3 / d, a b c d q^2 ; q^2\right)_{n-1}} \frac{\left(a c q, abq,bcq; q^2\right)_{n}}{(aq,b q, c q ; q)_{n}} \\
&-a~\frac{(1-d)(1-q/d)(1-bcq)(1-a^2bc)}{(1-abc)(1-a/d)(q-ad)}A_{-1}+q~\frac{(1-a/q)(1-a)(1-a b c q/ d)(1- a b c d)}{(1-abc)(1-a b c /q)(1-a/d)(q-ad)}\\
\times&\sum_{k=1}^{n}\left(1-a b c q^{3 k-1}\right)
\frac{(a b c/q,d, q/ d; q)_{k}\left(a c/q, ab/q,bcq; q^2\right)_{k}}
{\left(q^2,a b c q/ d, a b c d; q^2\right)_{k}(a/q,b q, c q ; q)_{k}}
 q^k,
\end{align*}
which is further reformulated  as
\begin{align}
&\sum_{k=0}^{n-1}\left(1-a b c q^{3k+1}\right) \frac{(abcq, d, q/ d ; q)_{k}
\left(b c q,acq,abq; q^2\right)_{k}}{\left(q^2, a b c q^3 / d, a b c d q^2 ; q^2\right)_k(aq,b q, c q ; q)_k} q^k \label{interestingxrma}\\
&=B_n(a)+C(a)\sum_{k=0}^{n}\left(1-a b c q^{3 k-1}\right)
\frac{(a b c/q,d, q/ d; q)_{k}\left(a c/q, ab/q,bcq; q^2\right)_{k}}
{\left(q^2,a b c q/ d, a b c d; q^2\right)_{k}(a/q,b q, c q ; q)_{k}}
 q^k,\nonumber
\end{align}
where
\begin{align}
C(a):=\frac{(1-a/q)(1-a)(1-a b c q/ d)(1- a b c d)}{(1-abc)(1-a b c /q)(1-a/d)(1-ad/q)}\label{coeefic};\end{align}
\begin{align}
&B_n(a):=\frac{a(1-q^{n+1}/a)(1-aq^n)(abc,d, q/ d ; q)_{n}\left(a c q, abq,bcq; q^2\right)_{n}}{(1-q^{2n})(1-abc)(1-a/d)(q-ad)\left(q^2,a b c q^3 / d, a b c d q^2 ; q^2\right)_{n-1}(aq,b q, c q ; q)_{n}} \nonumber \\
&-q~\frac{(1-a/q)(1-a)(1-a b c q/ d)(1- a b c d)}{(1-abc)(1-a/d)(q-ad)}-\frac{a(1-d)(1-q/d)(1-bcq)(1-a^2bc)}{(1 - a b c)(1 - a/d)(q- a d)}A_{-1}\nonumber.
\end{align}
At  this stage, we introduce the analytic function
\begin{align}F_n(a):=\sum_{k=0}^{n-1}\left(1-a b c q^{3k-1}\right)
\frac{(a b c/q,d, q/ d; q)_{k}\left(a c/q, ab/q,bcq; q^2\right)_{k}}
{\left(q^2,a b c q/ d, a b c d; q^2\right)_{k}(a/q,b q, c q ; q)_{k}}
 q^k,\end{align}
reformulating   \eqref{interestingxrma} in the form
 \begin{align}
F_n(aq^2)=C(a)F_{n+1}(a)+B_n(a).\label{analticrec}
\end{align}
By virtue of the basic relation \eqref{vvvv}, it is easily verified that
\begin{align*}q~\frac{(1-a/q)(1-a)(1-a b c q/ d)(1- a b c d)}{(1-abc)(1-a/d)(q-ad)}+a~\frac{(1-d)(1-q/d)(1-bcq)(1-a^2bc)}{(1 - a b c)(1 - a/d)(q- a d)}A_{-1}=0.\end{align*}
Thus $B_n(a)$  boils down to
\begin{align}
B_n(a)&=\frac{a(1-q^{n+1}/a)(1-aq^n)}{(1-q^{2n})(1-abc)(1-a/d)(q-ad)}\nonumber\\
&\times\frac{(abc,d, q/ d ; q)_{n}
}{\left(q^2,a b c q^3 / d, a b c d q^2 ; q^2\right)_{n-1}}\frac{\left(a c q, abq,bcq; q^2\right)_{n}}{(aq,b q, c q ; q)_{n}}\nonumber
\end{align}
and its limiting case $\lim_{n\to \infty}B_n(a)$, denoted by $B(a)$, is
\begin{align}
B(a)=\frac{a}{(1-a/d)(q-ad)}\frac{(abcq,d, q/ d ; q)_{\infty}
}{\left(q^2,a b c q^3 / d, a b c d q^2 ; q^2\right)_{\infty}}\frac{\left(a c q, abq,bcq; q^2\right)_{\infty}}{(aq,b q, c q ; q)_{\infty}}. \label{coeefib}
\end{align}
All that remains is to solve
\begin{align}
F_{\infty}(aq^2)=C(a)F_{\infty}(a)+B(a).\label{analticrec-infinite}
\end{align}
As a matter of fact, by iterating \eqref{analticrec-infinite} $n$ times, we have
\begin{align}
F_{\infty}(aq^{2(n+1)})&=F_{\infty}(a)\prod_{k=0}^nC(aq^{2k})+
\sum_{k=0}^nB(aq^{2(n-k)})\prod_{i=0}^{k-1}C(aq^{2(n-i)})\nonumber\\&=F_{\infty}(a)\prod_{k=0}^nC(aq^{2k})+
\sum_{k=0}^nB(aq^{2k})\prod_{i=k+1}^{n}C(aq^{2i}).\label{generalsol}
\end{align}
Since $F_{\infty}(a)$ is analytic in $a$, so we may take $n\to \infty$ and obtain
\begin{align*}
F_{\infty}(0)=F_{\infty}(a)\prod_{k=0}^\infty C(aq^{2k})+
\sum_{k=0}^\infty B(aq^{2k})\prod_{i=k+1}^\infty C(aq^{2i}).
\end{align*}
Finally we achieve
\begin{align}
F_{\infty}(a)&=\frac{F_{\infty}(0)}{\prod_{k=0}^\infty C(aq^{2k})}-
\sum_{k=0}^\infty B(aq^{2k})\frac{\prod_{i=k+1}^\infty C(aq^{2i})}{\prod_{i=0}^{\infty}C(aq^{2i})}\nonumber\\
&=F_{\infty}(0)\prod_{k=0}^\infty \frac{1}{C(aq^{2k})}-
\sum_{k=0}^\infty B(aq^{2k})\prod_{i=0}^{k}\frac{1}{C(aq^{2i})}.\label{yyyyyy}
\end{align}
In the meantime, we can calculate \begin{align*}
\prod_{i=0}^{k}\frac{1}{C(aq^{2i})}&=\prod_{i=0}^{k}\frac{(1-abcq^{2i})
(1-ab cq^{2i-1})(1-aq^{2i}/d)(1-aq^{2i-1}d)}{(1-aq^{2i-1})(1-aq^{2i})(1-ab c q^{2i+1}/ d)(1- a b c dq^{2i})}\\
&=\frac{(abc/q;q)_{2k+2}(a/d;q^2)_{k+1}(ad/q;q^2)_{k+1}}{(a/q;q)_{2k+2}(abcq/d;q^2)_{k+1}(abcd;q^2)_{k+1}}
\end{align*}
and the  limiting case
\begin{align*}
\prod_{k=0}^\infty \frac{1}{C(aq^{2k})}&=\frac{(abc/q;q)_{\infty}(a/d;q^2)_{\infty}(ad/q;q^2)_{\infty}}{(a/q;q)_{\infty}(abcq/d;q^2)_{\infty}(abcd;q^2)_{\infty}}.
\end{align*}
Substituting these expressions and $B(aq^{2k})$ into the recurrence \eqref{yyyyyy},  we finally obtain
\begin{align*}
F_{\infty}(a)&=\frac{(abc/q;q)_{\infty}(a/d;q^2)_{\infty}(ad/q;q^2)_{\infty}}{(a/q;q)_{\infty}(abcq/d;q^2)_{\infty}(abcd;q^2)_{\infty}}F_{\infty}(0)\\
&-a~
\frac{(abc/q,d, q/ d ; q)_{\infty}}{\left(q^2,a b c q/ d, a b c d; q^2\right)_{\infty}}\frac{\left(a c q, abq,bcq; q^2\right)_{\infty}}{(a/q,b q, c q ; q)_{\infty}}\sum_{k=0}^\infty\frac{(a/d,ad/q;q^2)_{k}}{\left(ac q, abq; q^2\right)_{k}}q^{2k-1}.
\end{align*}
So \eqref{gasperid-new} follows. It gives the complete proof of the theorem.
\qed

Two particular cases of Theorem \ref{thm3.2} deserve special mention.
\begin{tl}
\begin{align}
\sum_{k=0}^{\infty}
\frac{(d, q/ d; q)_{k}\left( ab/q; q^2\right)_{k}}
{\left(q^2; q^2\right)_{k}(a/q,b q; q)_{k}}
 q^k&=\frac{(a/d,ad/q;q^2)_{\infty}}{(a/q;q)_{\infty}}\sum_{k=0}^{\infty}
\frac{(d, q/ d; q)_{k}}
{\left(q^2; q^2\right)_{k}(b q; q)_{k}}
 q^k\label{gasperid-222}\\
&-a~
\frac{(d, q/ d ; q)_{\infty}}{\left(q^2; q^2\right)_{\infty}}\frac{\left(abq; q^2\right)_{\infty}}{(a/q,b q; q)_{\infty}}
\sum_{k=0}^\infty\frac{(a/d,ad/q;q^2)_{k}}{\left(abq; q^2\right)_{k}}q^{2k-1}.\nonumber
\end{align}
\end{tl}
\pf
  It suffices to set $c=0$ in \eqref{gasperid-new}.\qed
\begin{tl}
\begin{align}
\sum_{k=0}^{\infty}
\frac{(d, q/ d; q)_{k}}
{\left(q^2; q^2\right)_{k}(a/q; q)_{k}}
 q^k&=\frac{(a/d,ad/q;q^2)_{\infty}}{(a/q;q)_{\infty}}\sum_{k=0}^{\infty}
\frac{(d, q/ d; q)_{k}}
{\left(q^2; q^2\right)_{k}}
 q^k\label{gasperid-333}\\
&-a~
\frac{(d, q/ d ; q)_{\infty}}{\left(q^2; q^2\right)_{\infty}(a/q; q)_{\infty}}
\sum_{k=0}^\infty(a/d,ad/q;q^2)_{k}q^{2k-1}.\nonumber
\end{align}
\end{tl}
\pf It suffices to take $b=0$   in \eqref{gasperid-222}.\qed

Let us now consider the finite case of cubic $q$-series transformation, which was studied by Gasper  (cf. \cite[Eq. (5.22)]{Gasper89}).  Actually, by Theorem  \ref{type-i}  it is straightforward to verify
\begin{dl} Let define
\begin{align}
G_n(a,c):=\sum_{k=0}^{n-1}\frac{(aq^{4k}-1)q^{k}}{1-a} \frac{(a,b;q)_{k}}{(aq/c,cq/ab;q)_{k}}
 \frac{(c,{a^2b}/{c};q^3)_{k}}{(q^3,q^3a/b;q^3)_{k}}
            \frac{(q/b;q)_{2k}}{(ab;q)_{2k}}.\label{gabc}
\end{align}
Then
\begin{align}
G_n(a,c)&=\frac{\left(1-a b^2\right) (aq;q)_3}{(ab;q)_3(1-q^3a/b)}G_n(aq^3,cq^3)\label{pppppp-revised}\\
&+\frac{ab}{(1-a)(1-ab)}\frac{(a,b;q)_{n}}{(qa/c,qc/ab;q)_{n-1}} \frac{(q^3c,q^3a^2b/c;q^3)_{n-1}}{(q^3,q^3a^2b;q^3)_{n-1}}\frac{(q/b;q)_{2n}}{(qab;q)_{2n}}
      \frac{(q^3a^2b;q^3)_n}{(q^3a/b;q^3)_n}.\nonumber
\end{align}
\end{dl}
\pf It suffices to set
\begin{align*}
\left\{\begin{array}{l}
(p_1,p_2)\to (q,q),~~(p_3,p_4)\to(q^3,q^3)\\
(a_1^2,a_2^2,a_3^2,a_4^2)\to (aq,q/ab,cq^3,a^2bq^3/c);
\end{array}\right.
\left\{\begin{array}{l}
(q_1,q_4)\to (q,q^3),~~(q_2,q_3)\to(q^2,q^2)\\
(x_1^2,x_2^2,x_3^2,x_4^2)\to (b ,q/b,q^2/b,a^2bq^3).
\end{array}\right.
\end{align*}
Consequently, we easily find that
\begin{align*}
\left\{\begin{array}{lll}
K:=aq^4,&L:=q^4,&(L/p_1,L/p_2,L/p_3,L/p_4)=(q^3,q^3,q,q)\\
K_0:=aq^3,&L_0:=q^4,&(L_0/q_1,L_0/q_2,L_0/q_3,L_0/q_4)=(q^3,q^2,q^2,q)
\end{array}\right.
\end{align*}and
\begin{subequations}\label{kkkkkk2}
\begin{align}
\Gamma_{k-1}[\bar{a};\bar{p}]&=\frac{(c-1) (a^2 b-c)
}{b c}(1-abq^{2k})q^{5k};\\
\Gamma_{k}[\bar{x};\bar{q}]&=-\frac{\left(1-a b^2\right) \left(1-a q^{k+1}\right)\left(1-a q^{k+2}\right) }{b}q^{5k+3}
\end{align}as well  as
\begin{align*}
\prod_{i=1}^4(K/a_i^2;L/p_i)_k&=(qa/c,qc/ab;q)_k(q^3,q^3a^2b;q^3)_k;\\
\prod_{i=1}^4(K_0/x_i^2;L_0/q_i)_k&=(1/ab;q)_k(qab;q)_{2k}(q^3a/b;q^3)_k.
 \end{align*}
  \end{subequations}
Therefore,   \eqref{pppppp}  can be specificized to
 (note that $Kp_1=La_1^2=aq^5$)
\begin{align*}
&\sum_{k=0}^{n-1}\frac{(1-aq^{4k})q^{k}}{1-a}\frac{(a,1/ab;q)_{k}}{(qa/c,qc/ab;q)_{k}}\frac{(c,{a^2b}/{c};q^3)_{k}}{(q^3,q^3a^2b;q^3)_{k}}
            \frac{(b;q)_k}{(1/ab;q)_k}\frac{(q/b;q)_{2k}}{(ab;q)_{2k}}
      \frac{(q^3a^2b;q^3)_k}{(q^3a/b;q^3)_k}\nonumber\\
  &=\frac{(aq,q/ab;q)_{n-1}}{(qa/c,qc/ab;q)_{n-1}} \frac{(q^3c,q^3a^2b/c;q^3)_{n-1}}{(q^3,q^3a^2b;q^3)_{n-1}}
           \frac{(b;q)_n}{(1/ab;q)_n}\frac{(q/b;q)_{2n}}{(qab;q)_{2n}}
      \frac{(q^3a^2b;q^3)_n}{(q^3a/b;q^3)_n}\\
&+\frac{\left(1-a b^2\right) (aq;q)_3}{(1-q^3a/b)(ab;q)_3}\sum_{k=0}^{n-1}
\frac{(1-aq^{4k+3})q^{k}}{1-aq^3}\frac{(q^3a;q)_{k}}{(qa/c,qc/ab;q)_{k}}\frac{(q^3c,q^3a^2b/c;q^3)_{k}}{(q^3,q^3a^2b;q^3)_{k}} \nonumber\\
&\qquad\qquad\qquad\qquad\qquad\times
            \frac{(b;q)_k(q/b;q)_{2k}}{(q^3ab;q)_{2k}}
      \frac{(q^3a^2b;q^3)_k}{(q^6a/b;q^3)_k}\nonumber
.
\end{align*}
A further simplification gives rise to \eqref{pppppp-revised}.
\qed

By using \eqref{pppppp-revised}, we may recover a transformation found by  C. Y. Wang and J. N. Xu in \cite{xuwang} directly.
\begin{tl}[{\rm cf.  \cite[Sec. 2]{xuwang}}]Let $G_n(a,c)$ be given by \eqref{gabc}. Then\begin{align}
G_\infty(a,c)&=
G_\infty(q^3a,q^3c)\frac{(1-ab^2)(qa;q)_3}
{(1-q^3a/b)(ab;q)_3}\\
&+ab\frac{(aq,b,q/b;q)_\infty}{(aq/c,cq/ab,ab;q)_\infty}\frac{(q^3c,q^3a^2b/c;q^3)_\infty}{(q^3,aq^3/b;q^3)_\infty}.\nonumber
\end{align}
\end{tl}
\pf It suffices to put $n\to \infty $ on both sides of  \eqref{pppppp-revised}.
\qed

At the end of this subsection, we proceed to  the quartic transformation  \cite[Ex.3.33]{10} due to Gasper, which states
\begin{align*}
&\sum_{k=0}^{\infty} \frac{1-a^2 b^2 q^{5 k-2}}{1-a^2 b^2 / q^2} \frac{(a, b ; q)_k\left(a b / q, a b, a b q ; q^3\right)_k\left(a^2 b^2 / q^2 ; q^4\right)_k}{\left(a b^2 q^2, a^2 b q^2 ; q^4\right)_k\left(a b q, a b, a b / q ; q^2\right)_k(q ; q)_k} q^k\\
=&\frac{(a q, b ; q)_{\infty}\left(a^2 b^2 q^2 ; q^4\right)_{\infty}}{(q ; q)_{\infty}\left(a b q ; q^2\right)_{\infty}\left(b, a b^2 q^2, a^2 b q^2 ; q^4\right)_{\infty}}{ }_1 \phi_1\left[\begin{array}{c}
a \\
a q^4
\end{array}; q^4, b q^4\right].
\end{align*}
Inspired by this quartic transformation, we proceed to consider the following finite sum
\begin{align}
G_n(a,b):=\sum_{k=0}^{n-1}(1-a^2b^2q^{5k-2})\frac{(a, b ; q)_k\left(a^2 b^2/q^2 ; q^4\right)_{k}}{\left(a b^2 q^2, a^2 b q^2; q^4\right)_k(q; q)_k}
\frac{(abq,ab/q,ab;q^3)_k}{\left(a b q, a b, a b / q ; q^2\right)_k}q^k.\label{functiongab}
\end{align}
With the help of Theorem \ref{type-i}, we may establish the following new recurrence for $\{G_n(a,b)\}_{n\geq 0}$.
\begin{dl}Let $G_n(a,b)$ be given by \eqref{functiongab}. Then, for any integer $n\geq 0$, it holds
\begin{align}
G_n(a,b) &=G_{n-3}(aq^4,bq^4)\frac{ab(1-ab/q)q^4(aq, bq ; q)_{3}(-ab/q;q^2)_4}
{\left(1-a bq^7\right)\left(a b^2 q^2, a^2 b q^2 ; q^4\right)_{3}}\label{wantedresult}\\
&+
\frac{(aq, bq ; q)_{n-1}\left(a bq^2; q^4\right)_{n-1}\left(a^2 b^2/q^2; q^4\right)_{n}}{\left(a b^2 q^2, a^2 b q^2 ; q^4\right)_{n-1}(abq,q; q)_{n-1}}\frac{(abq;q)_n(ab/q^2,ab/q,ab;q^3)_n}
{\left(ab/q^2 ; q^4\right)_n\left(a b q, a b, a b / q ; q^2\right)_n}.\nonumber
\end{align}
\end{dl}
\pf  Given \eqref{pppppp}, we may choose
\begin{align*}
p_1=p_2=q, p_3=p_4=q^4,~~&
(a_1^2,a_2^2,a_3^2,a_4^2)=(aq,bq,abq^2,a^2b^2q^2);\\
 q_1=q,q_2=q_3=q_4=q^{3},~~&(x_1^2,x_2^2,x_3^2,x_4^2)
 =(abq,ab/q^2,ab/q,ab).
\end{align*}
Then it holds
\begin{align*}
K=a^2b^2q^3, L=q^5; \quad K_0=a^2b^2/q, L_0=q^5.
\end{align*}
It is easy to compute
\begin{align*}
\Gamma_{k-1}[\bar{a};\bar{p}]&:=-a^2b^2(1-a)(1-b)(1-abq^{3k-2})q^{6k-2};\\
\Gamma_{k}[\bar{x};\bar{q}]&:=a^3b^3(1-q^{k-2})(1-q^{k-1})(1-q^k)q^{6k}.
\end{align*}
In this case, we reduce \eqref{pppppp} to
\begin{align*}
&\sum_{k=0}^{n-1}(1-a^2b^2q^{5k-2})q^k
\frac{(a, b ; q)_k\left(a^2 b^2 q^2 ; q^4\right)_{k-1}}{\left(a b^2 q^2, a^2 b q^2 ; q^4\right)_k(q; q)_k}
\frac{({abq},ab/q,ab;q^3)_k}
{\left(a b q, a b, a b / q ; q^2\right)_k}\nonumber\\
  &=\frac{(aq, bq ; q)_{n-1}\left(a bq^2, a^2 b^2 q^2 ; q^4\right)_{n-1}}{\left(a b^2 q^2, a^2 b q^2 ; q^4\right)_{n-1}(abq,q; q)_{n-1}}\frac{(abq;q)_n(ab/q^2,ab/q,ab;q^3)_n}
{\left(ab/q^2 ; q^4\right)_n\left(a bq, a b, a b/q; q^2\right)_n}\\
+\frac{ab}{{1-ab/q^2}}&\sum_{k=0}^{n-1}(q^{k-2};q)_3(1-a^2b^2q^{5k-1})q^{k+1} \frac{({ab/q^2},ab/q,ab;q^3)_k}
{\left(a b q, a b, a b / q ; q^2\right)_{k+1}}\frac{(aq, bq ; q)_{k}\left(a^2 b^2 q^2 ; q^4\right)_{k}}{\left(a b^2 q^2, a^2 b q^2 ; q^4\right)_{k}(q; q)_{k}}.\nonumber
\end{align*}Observe that $(q^{k-2};q)_3=0, k=0,1,2$. After some simplifications, it becomes
\begin{align*}
&\sum_{k=0}^{n-1}(1-a^2b^2q^{5k-2})q^k
\frac{(a, b ; q)_k\left({a^2 b^2/q^2}; q^4\right)_{k}}{\left(a b^2 q^2, a^2 b q^2 ; q^4\right)_k(q; q)_k}
\frac{(abq,ab/q,ab;q^3)_k}
{\left(a b q, a b, a b / q ; q^2\right)_k}\nonumber\\
  &=\frac{(aq, bq ; q)_{n-1}\left(abq^2;q^4\right)_{n-1}\left(a^2 b^2/q^2 ; q^4\right)_{n}}{\left(a b^2 q^2, a^2 b q^2 ; q^4\right)_{n-1}(abq,q; q)_{n-1}}\frac{(abq;q)_n(ab/q^2,ab/q,ab;q^3)_n}
{\left(ab/q^2 ; q^4\right)_n\left(a b q, a b, a b / q ; q^2\right)_n}\\
&+ab\sum_{k=3}^{n-1}(1-a^2b^2q^{5k-1})q^{k+1}\frac{({abq};q^3)_{k-1}(ab/q,ab;q^3)_{k}}
{\left(a b q, a b, a b / q ; q^2\right)_{k+1}}
\frac{(aq, bq ; q)_{k}\left(a^2 b^2/q^2 ; q^4\right)_{k+1}}{\left(a b^2 q^2, a^2 b q^2 ; q^4\right)_{k}(q; q)_{k-3}}.\nonumber
\end{align*}
Then we find the recursive relation
\begin{align}
G_n(a,b) &=\frac{(aq, bq ; q)_{n-1}\left(abq^2;q^4\right)_{n-1}\left(a^2 b^2/q^2 ; q^4\right)_{n}}{\left(a b^2 q^2, a^2 b q^2 ; q^4\right)_{n-1}(abq,q; q)_{n-1}}\frac{(abq;q)_n(ab/q^2,ab/q,ab;q^3)_n}
{\left(ab/q^2 ; q^4\right)_n\left(a b q, a b, a b / q ; q^2\right)_n}\nonumber\\
&+L(a,b)\sum_{k=0}^{n-4}(1-a^2b^2q^{5k+14})q^{k}
\frac{(abq^7,abq^8,abq^9;q^3)_{k}}
{\left(a b q^9, a bq^8, a bq^7 ; q^2\right)_{k}}\frac{(aq^4, bq^4 ; q)_{k}\left(a^2 b^2 q^{14}; q^4\right)_{k}}{\left(a b^2 q^{14}, a^2 b q^{14} ; q^4\right)_{k}(q; q)_{k}},\label{newgab}
\end{align}
where $G_n(a,b)$ is given by \eqref{functiongab} and
 \begin{align*}L(a,b)&:=abq^4\frac{(abq;q^3)_2(ab/q,ab;q^3)_{3}}
{\left(a b q, a b, a b / q ; q^2\right)_{4}}\frac{(aq, bq ; q)_{3}\left(a^2 b^2/q^2 ; q^4\right)_{4}}{\left(a b^2 q^2, a^2 b q^2 ; q^4\right)_{3}}\\&=\frac{abq^4(aq, bq ; q)_{3}(1-ab/q)(-ab/q;q^2)_4}
{\left(1-a bq^7\right)\left(a b^2 q^2, a^2 b q^2 ; q^4\right)_{3}}.
\end{align*}
It is of interest to see that  the sum on the right side of \eqref{newgab} is just  $G_{n-3}(aq^4,bq^4)$.
Therefore, we have
\begin{align*}
G_n(a,b)&=\frac{(aq, bq ; q)_{n-1}\left(a bq^2;q^4\right)_{n-1}\left(a^2 b^2/q^2 ; q^4\right)_{n}}{\left(a b^2 q^2, a^2 b q^2 ; q^4\right)_{n-1}(abq,q; q)_{n-1}}\frac{(abq;q)_n(ab/q^2,ab/q,ab;q^3)_n}
{\left(ab/q^2 ; q^4\right)_n\left(a b q, a b, a b / q ; q^2\right)_n}\\
&+\frac{ab(1-ab/q)q^4(aq, bq ; q)_{3}(-ab/q;q^2)_4}
{\left(1-a bq^7\right)\left(a b^2 q^2, a^2 b q^2 ; q^4\right)_{3}}
G_{n-3}(aq^4,bq^4).\nonumber
\end{align*}
This gives the complete proof of the theorem.
\qed


\subsection{Some concrete transformations deduced from the second transformation}
In this  part, we willinvestigate some applications concerning Theorem \ref{type-ii}.
\begin{tl}\label{thm3.4}For any integer $n\geq 0$, it always holds
\begin{align}
&C_0\sum_{k=0}^{n-1}\frac{1-a b c q^{3 k+1}}{1-a}q^{k}\tau(k)
\frac{(d, q/ d ; q)_{k}
\left(b c q, a^2bc; q^2\right)_{k}}{\left(a b c q^3 / d, a b c d q^2 ; q^2\right)_k(a q; q)_k}\frac{ (a b cq; q)_k\left(a c q, abq; q^2\right)_k}{(abcq^2;q^3)_{k}(a;q)_{k}}\nonumber\\
&=\tau(n)(1-q^{n+1}/a)\frac{(d, q/ d ; q)_{n}
\left(b c q, a^2bc; q^2\right)_{n}}{\left(a b c q/ d, a b c d; q^2\right)_{n}(a; q)_{n}}
 \frac{(a b c; q)_n\left(a c q, abq; q^2\right)_n}{(abcq^2;q^3)_{n}(a;q)_{n}}\label{zzzz}
 \\
&-(1-q/a)-\sum_{k=1}^{n}\frac{\left(1-q^{2k}\right)\left(1-a^2bcq^{2k-2}\right)\left(1-bq^{k}\right)\left(1-cq^{k}\right)}{(1-abc/q)(1-ac/q)(1-ab/q)}\tau(k)\nonumber \\
&\qquad\qquad\qquad\times \frac{(d, q/ d ; q)_{k}
\left(b c q, a^2bc; q^2\right)_{k}}{\left(a b c q / d, a b c d; q^2\right)_{k}(a; q)_{k}}~\frac{(a b c/q; q)_k\left(a c/q, ab/q; q^2\right)_k}{(abcq^2;q^3)_{k}(a;q)_{k}},\nonumber
\end{align}
where\begin{align}C_0&:=\frac{(1-a/d)(q-ad)(1-abc)}{a(1-abcq/d)(1-abcd)};\label{xrma-1}\\
\tau(n)&:=\bigg(
\frac{-aq^{-(n-1)/2}}{(q-ab)(q-ac)} \bigg)^{n}.\label{xrma-2}\end{align}
\end{tl}
\pf
It suffices to make the parametric  replacement in \eqref{identityfinal-II} that
\begin{align*}
p_1=p_2=q,~~&p_3=p_4=q^2,~~(a_1^2,a_2^2,a_3^2,a_4^2)=(dq,q^2/d, bcq^3,a^2bcq^2).\\
 q_1=q_2=q,~~&q_3=q_4=q^2,~~(x_1^2,x_2^2,x_3^2,x_4^2)=(abc,q^2/a,acq, abq).
\end{align*}
In the sequel, it is easy to compute
\begin{align*}K:=abcq^4,&~~L=q^3,~~\Gamma_{k-1}[\bar{a};\bar{p}]=bc(1-a/d)(ad-q)(1-abcq^k)q^{4k+1};\\
K_0:=abcq^2,~~&L_0=q^3,~~\Gamma_{k}[\bar{x};\bar{q}]=bc(q-ab)(q-ac)(1-aq^k)q^{4k+2}.
\end{align*}
Accordingly, we easily find
\begin{align*}
A_k&:=\frac{(dq, q^2/ d ; q)_{k}
\left(b c q^3, a^2bcq^2 ; q^2\right)_{k}}{\left(a b c q^3 / d, a b c d q^2 ; q^2\right)_k(a q,q^2/a ; q)_k};\\
B_k&:=\frac{ (a b c,q^2/a; q)_k\left(a c q, abq; q^2\right)_k}{(abcq^2;q^3)_{k}(a;q)_{k}}\bigg(
\frac{-aq^{-(k-1)/2}}{(q-ab)(q-ac)} \bigg)^{k}.
\end{align*}
It is not difficult to compute the differences
\begin{align*}
\Delta A_k=\frac{(1-a/d)(q-ad)}{a}(1-abcq^k)(1-a b c q^{3 k+1}) \frac{(dq, q^2/ d ; q)_{k-1}
\left(b c q^3, a^2bcq^2 ; q^2\right)_{k-1}}{\left(a b c q^3 / d, a b c d q^2 ; q^2\right)_k(a q,q^2/a ; q)_k}
q^{k}\\
\end{align*}
and
\bnm
\nabla B_k
&=&\frac{(-1)^k}{(abcq^2;q^3)_{k+1}}\bigg(
\frac{a}{(q-ab)(q-ac)} \bigg)^{k+1}\frac{q^{-k(k+1)/2}}{(a;q)_{k+1}} \\
&\times&\left(1-q^{2k+2}\right)\left(1-a^2bcq^{2k}\right)\left(1-bq^{k+1}\right)\left(1-cq^{k+1}\right)\\
&\times& (a b c,q^2/a; q)_k\left(a c q, abq; q^2\right)_k\\
 &=&-\frac{(a b c,q^2/a; q)_k\left(a c q, abq; q^2\right)_k}{(abcq^2;q^3)_{k+1}(a;q)_{k+1}}\bigg(
\frac{-aq^{-k/2}}{(q-ab)(q-ac)} \bigg)^{k+1}\\
&\times&\left(1-q^{2k+2}\right)\left(1-a^2bcq^{2k}\right)\left(1-bq^{k+1}\right)\left(1-cq^{k+1}\right).
\enm
So we may specialize \eqref{abelparts} to
\begin{align}
&\frac{(1-a/d)(q-ad)}{a}\sum_{k=0}^{n-1}(1-abcq^k)\left(1-a b c q^{3 k+1}\right) \frac{(dq, q^2/ d ; q)_{k-1}
\left(b c q^3, a^2bcq^2 ; q^2\right)_{k-1}}{\left(a b c q^3 / d, a b c d q^2 ; q^2\right)_k(a q,q^2/a ; q)_k} q^k\nonumber
\\
&\times \frac{ (a b c,q^2/a; q)_k\left(a c q, abq; q^2\right)_k}{(abcq^2;q^3)_{k}(a;q)_{k}}\bigg(
\frac{-aq^{-(k-1)/2}}{(q-ab)(q-ac)} \bigg)^{k} \nonumber
\\
&=\frac{(dq, q^2/ d ; q)_{n-1}
\left(b c q^3, a^2bcq^2 ; q^2\right)_{n-1}}{\left(a b c q^3 / d, a b c d q^2 ; q^2\right)_{n-1}(a q,q^2/a ; q)_{n-1}} \frac{(a b c,q^2/a; q)_n\left(a c q, abq; q^2\right)_n}{(abcq^2;q^3)_{n}(a;q)_{n}}\bigg(
\frac{-aq^{-(n-1)/2}}{(q-ab)(q-ac)} \bigg)^{n}\label{xinrong-tomorrow1}\\
&-A_{-1}-\sum_{k=0}^{n-1}\left(1-q^{2k+2}\right)\left(1-a^2bcq^{2k}\right)\left(1-bq^{k+1}\right)\left(1-cq^{k+1}\right)\bigg(
\frac{-aq^{-k/2}}{(q-ab)(q-ac)} \bigg)^{k+1}\nonumber \\
&\qquad\qquad\qquad\times \frac{(dq, q^2/ d ; q)_{k}
\left(b c q^3, a^2bcq^2 ; q^2\right)_{k}}{\left(a b c q^3 / d, a b c d q^2 ; q^2\right)_k(a q,q^2/a ; q)_k}~\frac{(a b c,q^2/a; q)_k\left(a c q, abq; q^2\right)_k}{(abcq^2;q^3)_{k+1}(a;q)_{k+1}}.\nonumber
\end{align}
In view of  \eqref{vvvv}, we compute $A_{-1}$ and substitute it into the above identity.  The result is as follows.
\begin{align}
&\frac{(1-a/d)(q-ad)}{a}\sum_{k=0}^{n-1}(1-abcq^k)\left(1-a b c q^{3 k+1}\right) \frac{(dq, q^2/ d ; q)_{k-1}
\left(b c q^3, a^2bcq^2 ; q^2\right)_{k-1}}{\left(a b c q^3 / d, a b c d q^2 ; q^2\right)_k(a q; q)_k}\nonumber
\\
&\qquad\qquad\times \frac{ (a b c; q)_k\left(a c q, abq; q^2\right)_k}{(abcq^2;q^3)_{k}(a;q)_{k}}\bigg(
\frac{-aq^{-k/2+3/2}}{(q-ab)(q-ac)} \bigg)^{k}\label{zzzz-new}
\\
&=\frac{(dq, q^2/ d ; q)_{n-1}
\left(b c q^3, a^2bcq^2 ; q^2\right)_{n-1}}{\left(a b c q^3 / d, a b c d q^2 ; q^2\right)_{n-1}(a q,q^2/a ; q)_{n-1}}
 \frac{(a b c,q^2/a; q)_n\left(a c q, abq; q^2\right)_n}{(abcq^2;q^3)_{n}(a;q)_{n}}\bigg(
\frac{-aq^{-(n-1)/2}}{(q-ab)(q-ac)} \bigg)^{n}\nonumber\\
&\qquad\qquad-q~\frac{(1-1/a) \left(1-a/q\right) (1-a b c d) \left(1-a b c q/d\right)}{(1-d) \left(1-q/d\right) \left(1-a^2 b c\right) (1-b c q)}\nonumber \\
-&\sum_{k=0}^{n-1}\left(1-q^{2k+2}\right)\left(1-a^2bcq^{2k}\right)\left(1-bq^{k+1}\right)\left(1-cq^{k+1}\right)\bigg(
\frac{-aq^{-k/2}}{(q-ab)(q-ac)} \bigg)^{k+1}\nonumber \\
&\qquad\qquad \times\frac{(dq, q^2/ d ; q)_{k}
\left(b c q^3, a^2bcq^2 ; q^2\right)_{k}}{\left(a b c q^3 / d, a b c d q^2 ; q^2\right)_k(a q; q)_k}~\frac{(a b c; q)_k\left(a c q, abq; q^2\right)_k}{(abcq^2;q^3)_{k+1}(a;q)_{k+1}}.\nonumber
\end{align}
As the last step, by multiplying both sides of \eqref{zzzz-new} with
\[\frac{(1-d) \left(1-q/d\right) \left(1-a^2 b c\right) (1-b c q)}{(1-a)(1-a b c q/ d)(1- a b c d)},\]
we  achieve immediately
\begin{align}
&\frac{(1-a/d)(q-ad)(1-abc)}{a(1-a)(1-a b c q/ d)(1- a b c d)}\sum_{k=0}^{n-1}\left(1-a b c q^{3 k+1}\right)
 \frac{(d, q/ d ; q)_{k}
\left(b c q, a^2bc; q^2\right)_{k}}{\left(a b c q^3 / d, a b c d q^2 ; q^2\right)_k(a q; q)_k}\nonumber
\\
&\qquad\qquad\qquad\quad\times \frac{ (a b cq; q)_k\left(a c q, abq; q^2\right)_k}{(abcq^2;q^3)_{k}(a;q)_{k}}\bigg(
\frac{-aq^{-k/2+3/2}}{(q-ab)(q-ac)} \bigg)^{k}\label{zzzz-old}
\\
&=(1-q^{n+1}/a)\frac{(d, q/ d ; q)_{n}
\left(b c q, a^2bc; q^2\right)_{n}}{\left(a b c q/ d, a b c d; q^2\right)_{n}(a; q)_{n}}
 \frac{(a b c; q)_n\left(a c q, abq; q^2\right)_n}{(abcq^2;q^3)_{n}(a;q)_{n}}\bigg(
\frac{-aq^{-(n-1)/2}}{(q-ab)(q-ac)} \bigg)^{n}\nonumber \\
&-(1-q/a)-\sum_{k=0}^{n-1}\left(1-q^{2k+2}\right)\left(1-a^2bcq^{2k}\right)\left(1-bq^{k+1}\right)\left(1-cq^{k+1}\right)\bigg(
\frac{-aq^{-k/2}}{(q-ab)(q-ac)} \bigg)^{k+1}\nonumber \\
&\qquad\qquad\qquad\times \frac{(d, q/ d ; q)_{k+1}
\left(b c q, a^2bc; q^2\right)_{k+1}}{\left(a b c q / d, a b c d; q^2\right)_{k+1}(a; q)_{k+1}}~\frac{(a b c; q)_k\left(a c q, abq; q^2\right)_k}{(abcq^2;q^3)_{k+1}(a;q)_{k+1}}.\nonumber
\end{align}
 In view of \eqref{xrma-1} and  \eqref{xrma-2}, we  can finally express \eqref{zzzz-old} in the form of \eqref{zzzz}. The proof is finished.
\qed

Next are some special cases of Theorem \ref{thm3.4}.
\begin{lz}
\begin{align}
&\frac{(1-a/d)(q-ad)}{a(1-a)}\sum_{k=0}^{n-1} \frac{(d, q/ d ; q)_{k}
}{(a, a q; q)_k} \big(
-aq^{-k/2-1/2} \big)^{k}\label{xinrong-tomorrow333}
\\
&=\frac{(d, q/ d ; q)_{n}}{(a; q)_{n}^2}\big(
-aq^{-(n-1)/2-2}\big)^{n}(1-q^{n+1}/a)-\left(1-q/a\right)  \nonumber\\
&-\sum_{k=1}^{n} \frac{(d, q/ d ; q)_{k}
}{(a; q)_{k}^2}\left(1-q^{2k}\right)\big(
-aq^{-k/2-3/2}\big)^{k}.\nonumber
\end{align}
\end{lz}
\pf It is a direct consequence of  \eqref{zzzz} with $b=c=0$.\qed

\begin{lz}Define
\[T(k):=(-1)^kq^{-k(k-1)}\bigg(\frac{(q; q^2)_{k}}{(q^2; q^2)_{k}}\bigg)^2.\] Then
\begin{align}
\sum_{k=0}^{n-1}T(k) \frac{1-q}{1-q^{2k+2}}
=T(n)\frac{q^{-2n}-1}{1-q} -\sum_{k=1}^{n} T(k)\frac{q^{-2k}-q^{2k}}{1-q}.\label{xinrong-tomorrow444}
\end{align}
\end{lz}
\pf To show \eqref{xinrong-tomorrow444}, we first make the base replacement $q\to q^2$ in \eqref{xinrong-tomorrow333}. We have
\begin{align}
&\frac{(1-a/d)(q^2-ad)}{a(1-a)}\sum_{k=0}^{n-1} \frac{(d, q^2/ d ; q^2)_{k}
}{(a, a q^2; q^2)_k} \big(
-aq^{-k-1} \big)^{k}\label{combinal-id}
\\
&=\frac{(d, q^2/ d ; q^2)_{n}}{(a; q^2)_{n}^2}\big(
-aq^{-n-3}\big)^{n}(1-q^{2n+2}/a)-\left(1-q^2/a\right)  \nonumber\\
&-\sum_{k=1}^{n} \frac{(d, q^2/ d ; q^2)_{k}
}{(a; q^2)_{k}^2}\left(1-q^{4k}\right)\big(
-aq^{-k-3}\big)^{k}.\nonumber
\end{align}
And then let $a=q^2, d=q$ and express the resulted in view of $T(n)$. We finally have \eqref{xinrong-tomorrow444} at once.
\qed

\begin{remark} The limitation
$q\to{} _{-}1$ of \eqref{xinrong-tomorrow444}  offers a basic combinatorial identity:
\begin{align}
\sum_{k=0}^{n-1} (-16)^{n-k}\frac{1}{k+1}\binom{2 k}{k}^2 =4n\binom{2 n}{n}^2-8\displaystyle\sum_{k=1}^{n}(-16)^{n-k}\binom{2 k}{k}^2k.\label{xinrong-tomorrow888}
\end{align}
\end{remark}
Furthermore, appealing Theorem \ref{type-ii}, we can   establish
\begin{tl}
\begin{align}
&\sum_{k=0}^{n}\frac{(1-abq^{2k})(1-aq^{4k})}
{(1-ab)\left(1-aq^{k}\right)}\left(\frac{abq}{1-a b^2  }\right)^{k}q^{-k(k-1)/2}\label{type-II-concrete}\\
&\qquad\times \frac{(b;q)_k(q/b;q)_{2k}}{(aq^3;q^4)_k\left(a q^{2};q\right)_k}\frac{(1/ab;q)_{k}(c,a^2b/c;q^3)_{k}}{(aq/c,cq/ab;q)_k(q^3;q^3)_k}\nonumber\\
  &=\frac{(qa,q/ab;q)_{n}(q^3c,{q^3a^2b}/{c};q^3)_{n}}
{(q^3,q^3a^2b;q^3)_{n}(qa/c,qc/ab;q)_{n}}\frac{(b;q)_{n+1}(q/b;q)_{2n+2}(q^3a^2b;q^3)_{n+1}}{(aq^3;q^4)_{n+1}(aq,aq^2;q)_{n+1}}\bigg(\frac{ab}{1-ab^2}\bigg)^{n+1}
q^{-n(n+1)/2}\nonumber\\
&-\sum_{k=0}^{n}\frac{\left(1-aq^{3k+3}/b\right)\left(1-abq^{2k+1}\right)\left(1-abq^{2k+2}\right) \left(1-q^{k}/ab\right)}
{\left(1-aq^2\right)\left(1-aq^3\right)\left(1-aq^{k+1}\right)}\nonumber\\
&\qquad\times\left(\frac{ab}{1-a b^2 }\right)^{k+1}q^{-k(k+1)/2}
\frac{(b;q)_k(q/b;q)_{2k}}{(aq^7;q^4)_{k}(q^3;q^3)_k}\frac{(q/ab;q)_{k}(cq^3,a^2bq^3/c;q^3)_{k}}{(aq^3,aq/c,cq/ab;q)_k}.
\nonumber
\end{align}
\end{tl}
\pf In order to show \eqref{type-II-concrete}, we consider \begin{align*}
\left\{\begin{array}{l}
(p_1,p_2)\to (q,q),~~(p_3,p_4)\to(q^3,q^3)\\
(a_1^2,a_2^2,a_3^2,a_4^2)\to (aq,q/ab,cq^3,a^2bq^3/c);
\end{array}\right.
~\left\{\begin{array}{l}
(q_1,q_4)\to (q,q^3),~~(q_2,q_3)\to(q^2,q^2)\\
(x_1^2,x_2^2,x_3^2,x_4^2)\to (b ,q/b,q^2/b,a^2bq^3).
\end{array}\right.
\end{align*} It is easy to find \eqref{kkkkkk2} and then substitute it into
 \eqref{identityfinal-II}. After simplification, we have
\begin{align}
&\frac{(c-1) (a^2 b-c)
}{b c}\sum_{k=0}^{n-1}\frac{(1-abq^{2k})q^{5k}}
{(aq^3;q^4)_k}\frac{aq^{4k}-1}{aq^{4k}}~ \prod_{i=1}^4\frac{(a_i^2;p_i)_{k-1}(x_i^2;q_i)_k}{(aq^4/a_i^2;q^4/p_i)_k\nonumber
}\\
&\qquad\qquad\qquad\times\left(\frac{ab}{1-a b^2 }\right)^k \prod_{j=0}^{k-1}
\frac{q^{4j+3}}{\left(1-a q^{j+1}\right)\left(1-a q^{j+2}\right)q^{5j+3}}\nonumber\\
  &=H_n[a,b,c;q]\label{identityfinal-111-111-111}\\
&-\sum_{k=0}^{n-1}\frac{1}
{a^2q^{8k+6}(aq^3;q^4)_{k+1}}\left(\frac{ab}{1-a b^2 }\right)^{k+1} \prod_{j=0}^{k}
\frac{q^{4j+3}}{\left(1-a q^{j+1}\right)\left(1-a q^{j+2}\right)q^{5j+3}}\nonumber\\
&\qquad\qquad\qquad\times
\prod_{i=1}^4\frac{(a_i^2;p_i)_k(x_i^2;q_i)_k\left(aq^{4k+3}-x_i^2 q_i^k\right)}
{(aq^4/a_i^2;q^4/p_i)_k},\nonumber
\end{align}
Hereafter, for  brevity, we temporally write
\[H_n[a,b,c;q]:=\frac{1}{(aq^3;q^4)_n} \prod_{i=1}^4\frac{(a_i^2;p_i)_{n-1}(x_i^2;q_i)_n}
{(aq^4/a_i^2;q^4/p_i)_{n-1}}\prod_{j=0}^{n-1}
\frac{aq^{4j+3}}{\Gamma_j[\bar{x};\bar{q}]}.\]After some routine computation, we get
\begin{align*}
\frac{(c-1) (a^2 b-c)
}{ab c}&\sum_{k=0}^{n-1}\frac{(1-abq^{2k})(aq^{4k}-1)q^{k}}
{(aq^3;q^4)_k}~\left(\frac{ab}{1-a b^2  }\right)^k
\frac{q^{-k(k-1)/2}}{\left(a q,a q^{2};q\right)_k} \\
&\times\frac{(aq,q/ab;q)_{k-1}(cq^3,a^2bq^3/c;q^3)_{k-1}}{(aq/c,cq/ab;q)_k(q^3,a^2bq^3;q^3)_k}(b;q)_k(q/b;q)_{2k}(a^2bq^3;q^3)_k\nonumber\\
  =H_n[a,b,c;q]-&\sum_{k=0}^{n-1}\frac{1}
{a^2q^{8k+6}(aq^3;q^4)_{k+1}}\left(\frac{ab}{1-a b^2 }\right)^{k+1}
\frac{q^{-k(k+1)/2}}{\left(a q,a q^{2};q\right)_{k+1}}\\
&\times\frac{(aq,q/ab;q)_{k}(cq^3,a^2bq^3/c;q^3)_{k}}{(aq/c,cq/ab;q)_k(q^3,a^2bq^3;q^3)_k}(b;q)_k(q/b;q)_{2k}(a^2bq^3;q^3)_k\\
&\times\left(aq^{4k+3}-b q^k\right)\left(aq^{4k+3}- q^{2k+1}/b\right)\left(aq^{4k+3}-q^{2k+2}/b\right)\left(aq^{4k+3}-a^2b q^{3k+3}\right).
\nonumber
\end{align*}
It is equivalent to
\begin{align*}
&\frac{1}{ab}\sum_{k=0}^{n-1}\frac{(1-abq^{2k})(aq^{4k}-1)q^{k}}
{(aq^3;q^4)_k}\frac{(b;q)_k(q/b;q)_{2k}}{\left(a q,a q^{2};q\right)_k} \left(\frac{ab}{1-a b^2  }\right)^{k}q^{-k(k-1)/2}\\
&\qquad\times\frac{(aq,q/ab;q)_{k-1}(c,a^2b/c;q^3)_{k}}{(aq/c,cq/ab;q)_k(q^3;q^3)_k}\nonumber\\
  &= H_n[a,b,c;q]\\
&-\sum_{k=0}^{n-1}\frac{\left(1-aq^{3k+3}/b\right)\left(1-abq^{2k+1}\right)\left(1-abq^{2k+2}\right) \left(1-q^{k}/ab\right)}
{(aq^3;q^4)_{k+1}}
\frac{(b;q)_k(q/b;q)_{2k}}{\left(a q,a q^{2};q\right)_{k+1}}\\
&\qquad\times\frac{(aq,q/ab;q)_{k}(cq^3,a^2bq^3/c;q^3)_{k}}{(aq/c,cq/ab;q)_k(q^3;q^3)_k}\left(\frac{ab}{1-a b^2 }\right)^{k+1}q^{-k(k+1)/2}.
\nonumber
\end{align*}
After some arrangement, we get
\begin{align*}
&\sum_{k=0}^{n-1}\frac{(1-abq^{2k})(1-aq^{4k})}
{(1-ab)\left(1-aq^{k}\right)}\left(\frac{abq}{1-a b^2  }\right)^{k}q^{-k(k-1)/2}\\
&\qquad\times \frac{(b;q)_k(q/b;q)_{2k}}{(aq^3;q^4)_k\left(a q^{2};q\right)_k}\frac{(1/ab;q)_{k}(c,a^2b/c;q^3)_{k}}{(aq/c,cq/ab;q)_k(q^3;q^3)_k}\nonumber\\
  &= H_n[a,b,c;q]\\
&-\sum_{k=0}^{n-1}\frac{\left(1-aq^{3k+3}/b\right)\left(1-abq^{2k+1}\right)\left(1-abq^{2k+2}\right) \left(1-q^{k}/ab\right)}
{\left(1-aq^2\right)\left(1-aq^3\right)\left(1-aq^{k+1}\right)}\\
&\qquad\times\left(\frac{ab}{1-a b^2 }\right)^{k+1}q^{-k(k+1)/2}
\frac{(b;q)_k(q/b;q)_{2k}}{(aq^7;q^4)_{k}(q^3;q^3)_k}\frac{(q/ab;q)_{k}(cq^3,a^2bq^3/c;q^3)_{k}}{(aq^3,aq/c,cq/ab;q)_k}.
\nonumber
\end{align*} The proof is finished.
\qed

\begin{lz}
\begin{align}
&\sum_{k=0}^{\infty}\frac{(1-abq^{2k})(1-aq^{4k})}
{(1-ab)\left(1-aq^{k}\right)}\left(\frac{a^2b^2}{1-a b^2  }\right)^{k}q^{k(k-1)/2} \frac{(b;q)_k(q/b;q)_{2k}}{(aq^3;q^4)_k\left(a q^{2};q\right)_k}\frac{(1/ab;q)_{k}}{(q^3;q^3)_k}\nonumber\\
=&-\frac{1}{ab}\sum_{k=0}^{\infty}\frac{\left(1-aq^{3k+3}/b\right)\left(1-abq^{2k+1}\right)\left(1-abq^{2k+2}\right) \left(1-q^k/ab\right)}
{\left(1-aq^2\right)\left(1-aq^3\right)\left(1-aq^{k+1}\right)}\label{3.47}\\
&\qquad\qquad\times\left(\frac{a^2b^2}{1-a b^2 }\right)^{k+1}q^{k(k+1)/2}
\frac{(b;q)_k(q/b;q)_{2k}}{(aq^7;q^4)_{k}(q^3;q^3)_k}\frac{(q/ab;q)_{k}}{(aq^3;q)_k}.
\nonumber
\end{align}
\end{lz}
\pf To show \eqref{3.47}, we need only  to put $c\to 0$ and then take $n\to\infty$ in \eqref{type-II-concrete} . Note that
\[\lim_{c\to 0}\frac{(A/c;q^3)_{k}}{(B/c;q)_k}=\bigg(\frac{A}{B}\bigg)^kq^{k(k-1)}.\]
Then the conclusion follows.
\qed

\subsection{Some specific transformations deduced from the third transformation}
In what follows, we will  present a special instance of Theorem \ref{type-iii}.
\begin{tl}\label{thm3.5}For any integer $n\geq 0$, it always holds
\begin{align}&\frac{(yq,wq,uq,yz^2wuq/x^2;q)_{n+1}}{
(q u w z/x,q u y z/x,q w y z/x,q x/z;q)_{n+1}}\bigg(1+\frac{a}{c}\frac{\mathcal{D} (bc/a,cd/,ce/a,b c d e/a)}{\mathcal{D}(c d e/a,b c e/a,b c d/a,a/c)}\nonumber\\
&
 \quad\times\sum_{k=0}^n\frac{1-b c d e p^{2 k}/a}{1-b c d e/a}p^k \frac{(b,d,e,bc^2de/a^2;p)_{k}}{(c d e p/a,b c e p/a,b c d p/a,a p/c;p)_k}\nonumber\\ &\qquad\qquad\times\frac{(x q^{-1-n}/u w z,x q^{-1-n}/u y z,x q^{-1-n}/ w y z,z q^{-1-n}/x;q)_{k}}{(q^{-1-n}/ y,q^{-1-n}/ w,q^{-1-n}/ u,q^{-1-n} x^2/ u w y z^2;q)_{k}}\bigg)
\nonumber\\
&=\frac{(b,d,e,bc^2de/a^2;p)_{n+1}}{(c d e/a,b c e /a,b c d/a,a/c;p)_{n+1}}\bigg(1+\frac{xq}{z} \frac{\mathcal{D} (yz/x,zu/x,zw/x,q^{2} u w y z/x)}{\mathcal{D}(q u w z/x,q u y z/x,q w y z/x,q x/z))}\nonumber\\
&\quad\times\sum_{k=0}^n\frac{1-q^{2+2 k} u w y z/x}{1-q^{2} u w y z/x}q^k\frac{(yq,wq,uq,yz^2wuq/x^2;q)_k}
{(q^2 u w z/x,q^2 u y z/x,q^2 w y z/x,q^2 x/z;q)_{k}}
\label{thirdadded}\\
&\qquad\qquad\times  \frac{(a p^{-n}/c d e,a p^{-n}/b c e ,a p^{-n}/b c d,c p^{-n}/a ;p)_{k}}{(p^{-n}/b ,p^{-n}/d ,p^{-n}/e,a^2 p^{-n}/b c^2 d e;p)_{k}}\bigg).\nonumber
\end{align}
\end{tl}
\pf To establish \eqref{thirdadded}, we specialize \eqref{alpha-beta-new} with the same parameters as those used by Theorem \ref{type-i}, namely
\begin{align*}
p_i=p, ~~(1\leq i\leq 4),~&
(a_1^2,a_2^2,a_3^2,a_4^2)=(bp,dp,ep,bc^2dep/a^2);\\
 q_i=q,~~(1\leq i\leq 4),~&(x_1^2,x_2^2,x_3^2,x_4^2)=(yq,wq,uq,yz^2wuq/x^2).
\end{align*}
These yield
\begin{align*}&K:=bcdep^2/a,~~L=p^2,\\
&\Gamma_{k-1}[\bar{a};\bar{p}]=-bde(1-bc/a)(1-cd/a)(1-ce/a)p^{3k};\\
&K_0:=yzwuq^2/x,~~L_0=q^2,\\
&\Gamma_{k}[\bar{x};\bar{q}]=-ywu(1-yz/x)(1-zw/x)(1-zu/x)q^{3k+3}.
\end{align*}
 A bit long simplification gives rise to
\begin{align*}\frac{(yq,wq,uq,yz^2wuq/x^2;q)_{n+1}}{
(q u w z/x,q u y z/x,q w y z/x,q x/z;q)_{n+1}}\bigg(1\\
+a/c (1-bc/a) (1-cd/a) (1-ce/a) \sum_{k=1}^n(1-b c d e p^{2 k}/a)p^k\nonumber\\ \times \frac{(bp,dp,ep,bc^2dep/a^2;p)_{k-1}}{(c d e p/a,b c e p/a,b c d p/a,a p/c;p)_k}\frac{(x q^{-1-n}/u w z,x q^{-1-n}/(u y z,x q^{-1-n}/ w y z,z q^{-1-n}/x;q)_{k}}{(q^{-1-n}/y,q^{-1-n}/w,q^{-1-n}/u,q^{-1-n} x^2/u w y z^2;q)_{k}}\bigg)
\nonumber\\
=\frac{(bp,dp,ep,bc^2dep/a^2;p)_{n}}{(c d e p/a,b c e p/a,b c d p/a,a p/c;p)_{n}}\bigg(1\\+xq/z (1-yz/x) (1-zu/x) (1-zw/x)\sum_{k=0}^n(1-q^{2+2 k} u w y z/x)q^k
\\
\qquad\times \frac{(yq,wq,uq,yz^2wuq/x^2;q)_k}
{(q u w z/x,q u y z/x,q w y z/x,q x/z;q)_{k+1}} \frac{(a p^{-n}/c d e,a p^{-n}/b c e,a p^{-n}/b c d,c p^{-n}/a ;p)_{k}}{(p^{-n}/b ,p^{-n}/d ,p^{-n}/e,a^2 p^{-n}/b c^2 d e;p)_{k}}\bigg),\nonumber
\end{align*}
which in turn can be reformulate as the form
\begin{align*}\frac{(yq,wq,uq,yz^2wuq/x^2;q)_{n+1}}{
(q u w z/x,q u y z/x,q w y z/x,q x/z;q)_{n+1}}\bigg(\frac{(1-a/c) (1-b c d/a) (1-b c e/a) (1-c d e/a)}{(1-b)  (1-d) (1-e) (1-b c^2 d e/a^2)}\\
+\frac{a}{c}\frac{ (1-bc/a) (1-cd/a) (1-ce/a)(1-b c d e/a)}{(1-b)(1-d)(1-e)(1-bc^2de/a^2)}
\sum_{k=0}^n\frac{1-b c d e p^{2 k}/a}{1-b c d e/a}p^k\nonumber\\ \times \frac{(b,d,e,bc^2de/a^2;p)_{k}}{(c d e p/a,b c e p/a,b c d p/a,a p/c;p)_k}\frac{(x q^{-1-n}/u w z,x q^{-1-n}/u y z,x q^{-1-n}/w y z,z q^{-1-n}/x;q)_{k}}{(q^{-1-n}/y,q^{-1-n}/w,q^{-1-n}/u,(q^{-1-n} x^2)/u w y z^2;q)_{k}}\bigg)
\nonumber\\
=\frac{(bp,dp,ep,bc^2dep/a^2;p)_{n}}{(c d e p/a,b c e p/a,b c d p/a,a p/c;p)_{n}}\bigg(1\\+\frac{xq}{z}\frac{ (1-yz/x) (1-zu/x) (1-zw/x)(1- u w y zq^{2}/x)}{(q u w z/x,q u y z/x,q w y z/x,q x/z)}\sum_{k=0}^n\frac{1- u w y zq^{2+2 k}/x}{1- u w y zq^{2}/x}q^k
\\
\qquad\times \frac{(yq,wq,uq,yz^2wuq/x^2;q)_k}
{(u w zq^2/x,u y zq^2/x,w y zq^2/x,xq^2/z;q)_{k}} \frac{(a p^{-n}/c d e,a p^{-n}/b c e,a p^{-n}/b c d ,c p^{-n}/a ;p)_{k}}{(p^{-n}/b ,p^{-n}/d ,p^{-n}/e,a^2 p^{-n}/b c^2 d e;p)_{k}}\bigg).\nonumber
\end{align*}
By dividing both sides by
\[\frac{(1-a/c) (1-b c d/a) (1-b c e/a) (1-c d e/a)}{(1-b)  (1-d) (1-e) (1-b c^2 d e/a^2)},\]
then we get \eqref{thirdadded} as desired.
\qed

It is of interest to see that if $x=zw$ in \eqref{thirdadded}, then it follows Gasper's bibasic  summation formula:
\begin{lz}
\begin{align}
&\sum_{k=0}^n\frac{1-b c d e p^{2 k}/a}{1-b c d e/a}\frac{(b,d,e,bc^2de/a^2;p)_{k}}{(a p/c,c d e p/a,b c e p/a,b c d p/a;p)_k}p^k \\
&=\frac{(1-c/a)(1-c d e/a)(1-b c e/a)(1-b c d/a)}{ (1-bc/a) (1-cd/a) (1-ce/a)(1-b c d e/a)} \bigg(1-\frac{(b,d,e,bc^2de/a^2;p)_{n+1}}{(a/c,c d e/a,b c e /a,
b c d/a;p)_{n+1}}\bigg).
\nonumber
\end{align}
\end{lz}
Furthermore, we can derive a transformation for ${}_{10}\phi_9$ series which is  evidently different from Bailey's transformation  for four very-well-poised ${}_{10}\phi_9$ series \cite[Eq.(2.12.9)]{10}.
\begin{lz}[{\rm\cite[Eq.(2.11)]{Gasper89}}]For integer $n\geq 0$, it holds
\begin{align}&{}_{10}\phi_{9}\left[\begin{matrix}bde,q\sqrt{bde},-q\sqrt{bde},b,d,e,q^{-2-n}/u w, q^{-2-n}/u y, q^{-2-n}/ w y,q^{-n}\\ \sqrt{bde},-\sqrt{bde}, d e q,b  e q,b  d q,q^{-1-n}/ y,q^{-1-n}/ w,q^{-1-n}/ u,q^{-3-n}/ u w y\end{matrix};q,q\right]
\nonumber\\
&=\frac{(bq,dq,eq,bdeq,u w q^3,u y q^3,w y q^3;q)_{n}}{( d eq,b  eq ,b dq,yq^{2},wq^{2},uq^{2},ywuq^4;q)_{n}} \label{succeed}\\
&\times{}_{10}\phi_{9}\left[\begin{matrix}ywuq^3,q^2\sqrt{ywuq},-q^2\sqrt{ywuq},yq,wq,uq,q^{-n}/ d e, q^{-n}/b e , q^{-n}/b d, q^{-n}\\ q\sqrt{ywuq},-q\sqrt{ywuq}, q^3 u w,q^3 u y,q^3 w y,q^{-n}/b ,q^{-n}/d ,q^{-n}/e,q^{-n}/b  d e\end{matrix};q,q\right].\nonumber
\end{align}
\end{lz}
\pf It suffice to take $p=q$ and then multiply both sides of \eqref{thirdadded} with $(1-a/c)(1-qx/z).$ We have
\begin{align*}&\frac{(yq,wq,uq,yz^2wuq/x^2;q)_{n+1}}{
(q u w z/x,q u y z/x,q w y z/x;q)_{n+1}(q^2x/z;q)_{n}}\bigg(1-a/c+\frac{a}{c}\frac{\mathcal{D} (bc/a,cd/,ce/a,b c d e/a)}{\mathcal{D}(c d e/a,b c e/a,b c d/a)}\nonumber\\
&
 \quad\times\sum_{k=0}^n\frac{1-b c d e q^{2 k}/a}{1-b c d e/a}q^k \frac{(b,d,e,bc^2de/a^2;q)_{k}}{(c d e q/a,b c e q/a,b c d q/a,a q/c;q)_k}\nonumber\\ &\qquad\qquad\times\frac{(x q^{-1-n}/u w z,x q^{-1-n}/u y z,x q^{-1-n}/ w y z,z q^{-1-n}/x;q)_{k}}{(q^{-1-n}/ y,q^{-1-n}/ w,q^{-1-n}/ u,q^{-1-n} x^2/ u w y z^2;q)_{k}}\bigg)
\nonumber\\
&=\frac{(b,d,e,bc^2de/a^2;q)_{n+1}}{(c d e/a,b c e /a,b c d/a;q)_{n+1}(aq/c;q)_{n}}\bigg(1-xq/z+\frac{xq}{z} \frac{\mathcal{D} (yz/x,zu/x,zw/x,q^{2} u w y z/x)}{\mathcal{D}(q u w z/x,q u y z/x,q w y z/x))}\nonumber\\
&\quad\times\sum_{k=0}^n\frac{1-q^{2+2 k} u w y z/x}{1-q^{2} u w y z/x}q^k\frac{(yq,wq,uq,yz^2wuq/x^2;q)_k}
{(q^2 u w z/x,q^2 u y z/x,q^2 w y z/x,q^2 x/z;q)_{k}}\nonumber\\
&\qquad\qquad\times  \frac{(a q^{-n}/c d e,a q^{-n}/b c e ,a q^{-n}/b c d,c q^{-n}/a ;q)_{k}}{(q^{-n}/b ,q^{-n}/d ,q^{-n}/e,a^2 q^{-n}/b c^2 d e;q)_{k}}\bigg).\nonumber
\end{align*}
And then let $c=a$ and $z=xq$ in the last transformation, leading to
\begin{align*}&\sum_{k=0}^n\frac{1-bd e q^{2 k}}{1-b  d e}q^k \frac{(b,d,e,bde;q)_{k}}{( d e q,b  e q,b  d q,q;q)_k}\times\frac{(q^{-2-n}/u w, q^{-2-n}/u y, q^{-2-n}/ w y,q^{-n};q)_{k}}{(q^{-1-n}/ y,q^{-1-n}/ w,q^{-1-n}/ u,q^{-3-n}/ u w y;q)_{k}}
\nonumber\\
&=\frac{(bq,dq,eq,bdeq,u w q^3,u y q^3,w y q^3;q)_{n}}{( d eq,b  eq ,b dq,yq^{2},wq^{2},uq^{2},ywuq^4;q)_{n}} \nonumber\\
&\quad\times\sum_{k=0}^n\frac{1-q^{3+2 k} u w y}{1-q^{3} u w y}q^k\frac{(yq,wq,uq,ywuq^3;q)_k}
{(q^3 u w,q^3 u y,q^3 w y ,q;q)_{k}}\times  \frac{(q^{-n}/ d e, q^{-n}/b e , q^{-n}/b d, q^{-n} ;q)_{k}}{(q^{-n}/b ,q^{-n}/d ,q^{-n}/e,q^{-n}/b  d e;q)_{k}}.\nonumber
\end{align*}
It is just \eqref{succeed}.
\qed
\section{A unified treatment for  Gasper-Rahman's multibasic transformations}
In order to   better understand Gasper-Rahman's multibasic transformations of $q$-series   from a broader
perspective,  we think that it is necessary to explore further the special and more explicit form of Theorem \ref{type-i} with $p_i=p^{2r_i}$ and $q_i=q^{2s_i}$, $1\leq i\leq 4$.
To that end,  we first  deduce   a general multibasic transformation from Theorem \ref{type-i} .
\begin{dl}\label{theorem41}Let $\mathcal{D}$ be given by \eqref{charafunchi}. Then, for any integer $n$ and nonnegative real numbers $r_i,s_i$,$1\leq i\leq 4$, $r=r_1+r_2+r_3+r_4$, $s=s_1+s_2+s_3+s_4$, it holds
\begin{multline}\begin{split}
 \frac{a}{c}\sum_{k=0}^{n-1}\mathcal{D}\bigg(\frac{ce}{a}p^{(r_3+r_4-r_1-r_2)(k-1)},\frac{bc}{a}p^{(r_1+r_4-r_2-r_3)(k-1)},\frac{cd}{a}p^{(r_2+r_4-r_1-r_3)(k-1)},\frac{b c d e}{a} p^{r(k-1)+2} \bigg)\\
\times p^{(r-2r_4)(k-1)+1}
\frac{(bp;p^{2r_1})_{k-1}(dp;p^{2r_2})_{k-1}(ep;p^{2r_3})_{k-1}
(bc^2dep/a^2;p^{2r_4})_{k-1}}
{(cdep/a;p^{r-2r_1})_{k}(bcep/a;p^{r-2r_2})_{k}(bcdp/a;p^{r-2r_3})_{k}
(ap/c;p^{r-2r_4})_{k}}\\
\times\frac{(yq;q^{2s_1})_{k}(wq;q^{2s_2})_{k}(uq;q^{2s_3})_{k}
(yz^2wuq/x^2;q^{2s_4})_{k}}
{(zwuq/x;q^{s-2s_1})_{k}(yzuq/x;q^{s-2s_2})_{k}
(yzwq/x;q^{s-2s_3})_{k}(xq/z;q^{s-2s_4})_{k}}\\
=
\frac{(bp;p^{2r_1})_{n-1}(dp;p^{2r_2})_{n-1}(ep;p^{2r_3})_{n-1}
(bc^2dep/a^2;p^{2r_4})_{n-1}}
{(cdep/a;p^{r-2r_1})_{n-1}(bcep/a;p^{r-2r_2})_{n-1}(bcdp/a;p^{r-2r_3})_{n-1}
(ap/c;p^{r-2r_4})_{n-1}}\\
\times\frac{(yq;q^{2s_1})_{n}(wq;q^{2s_2})_{n}(uq;q^{2s_3})_{n}
(yz^2wuq/x^2;q^{2s_4})_{n}}
{(zwuq/x;q^{s-2s_1})_{n}(yzuq/x;q^{s-2s_2})_{n}
(yzwq/x;q^{s-2s_3})_{n}(xq/z;q^{s-2s_4})_{n}}\\
-\frac
{(p^{r-2r_1}-cdep/a)(p^{r-2r_2}-bcep/a)(p^{r-2r_3}-bcdp/a)(p^{r-2r_4}-ap/c)}{(p^{2r_1}-bp)(p^{2r_2}-dp)(p^{2r_3}-ep)(p^{2r_4}-bc^2dep/a^2)}\\
-\frac{x}{z}\sum_{k=0}^{n-1}\mathcal{D}\bigg(\frac{yz}{x}q^{(s_1+s_4-s_2-s_3)k},\frac{zu}{x}q^{(s_3+s_4-s_1-s_2)k},\frac{zw}{x}q^{(s_2+s_4-s_1-s_3)k},\frac{y z w u}{x} q^{s k+2}\big)\\\label{pppppp-generalid}\times q^{(s-2s_4)k+1}\frac{(yq;q^{2s_1})_{k}(wq;q^{2s_2})_{k}(uq;q^{2s_3})_{k}
(yz^2wuq/x^2;q^{2s_4})_{k}}
{(zwuq/x;q^{s-2s_1})_{k+1}(yzuq/x;q^{s-2s_2})_{k+1}
(yzwq/x;q^{s-2s_3})_{k+1}(xq/z;q^{s-2s_4})_{k+1}}\\
\times
\frac{(bp;p^{2r_1})_{k}(dp;p^{2r_2})_{k}(ep;p^{2r_3})_{k}
(bc^2dep/a^2;p^{2r_4})_{k}}
{(cdep/a;p^{r-2r_1})_{k}(bcep/a;p^{r-2r_2})_{k}(bcdp/a;p^{r-2r_3})_{k}
(ap/c;p^{r-2r_4})_{k}}.\end{split}
 \end{multline}
 \end{dl}
\pf It suffices to set in \eqref{pppppp} that
\begin{align*}
&p_i=p^{2r_i},  q_i=q^{2s_i}, ~1\leq i\leq j\leq 4,\\
&(a_1^2,a_2^2,a_3^2,a_4^2)=(bp,dp,ep,bc^2dep/a^2),~~
(x_1^2,x_2^2,x_3^2,x_4^2)=(yq,wq,uq,yz^2wuq/x^2).
\end{align*}
Subsequently, we get
\begin{align*}K:=bcdep^2/a,&~~L=p^r, r=r_1+r_2+r_3+r_4;\\
\Gamma_{k-1}[\bar{a};\bar{p}]&=-bdep^{2(r-r_4)(k-1)+3}\big(1-\frac{ce}{a}p^{(r_3+r_4-r_1-r_2)(k-1)}\big)\\
&\times(1-\frac{bc}{a}p^{(r_1+r_4-r_2-r_3)(k-1)}\big)\big(1-\frac{cd}{a}p^{(r_2+r_4-r_1-r_3)(k-1)}\big);\\
K_0:=yzwuq^2/x,~~&L_0=q^s, s=s_1+s_2+s_3+s_4;\\
\Gamma_{k}[\bar{x};\bar{q}]&=-ywuq^{2(s-s_4)k+3}\big(1-\frac{yz}{x}q^{(s_1+s_4-s_2-s_3)k}\big)\\
&\times\big(1-\frac{zu}{x}q^{(s_3+s_4-s_1-s_2)k}\big)\big(1-\frac{zw}{x}q^{(s_2+s_4-s_1-s_3)k}\big).
\end{align*}
Note that in this case, it is easy to check
\begin{align}
A_{-1}=\frac
{(p^{r-2r_1}-cdep/a)(p^{r-2r_2}-bcep/a)(p^{r-2r_3}-bcdp/a)(p^{r-2r_4}-ap/c)}{(p^{2r_1}-bp)(p^{2r_2}-dp)(p^{2r_3}-ep)(p^{2r_4}-bc^2dep/a^2)}.\label{A-1}
\end{align}
A direct substitution of these expressions into \eqref{pppppp} gives \eqref{pppppp-generalid}. The  detailed computation is left to the reader.
\qed

Following up on Gasper and Rahman \cite{Gasper90}, if $(a;q^m)_k$ occurs to a $q$-series transformation and $m$ is the maximum positive integer, then it is said to be of degree $m$. In this sense, when $m=2,3,4,5$, they call the corresponding  transformations  the quadratic, cubic, quartic, and quintic transformations. To generalize this  idea, we need to introduce

\begin{dy}[$(R,S)$-type transformation with degree $2m$] The transformation  \eqref{pppppp-generalid} is said to be $(R,S)$-type with degree $2m$, provided that $m=\max\{r_i,s_i|1\leq i\leq 4\}$ and $R,S$ are the cardinalities of sets $\{r_i|1\leq i\leq 4\}$ and  $\{s_i|1\leq i\leq 4\}$, respectively.
\end{dy}

From this definition, we see that \eqref{pppppp-generalid} is of the $(4,4)$-type with degree $2m$. In what follows, we  will detail  some special transformations with an arbitrary degree $2m$ for further use.

\begin{tl}[$(1,1)$-type transformation with degree $2m$] \label{degree-2s2-1-new}Let $\mathcal{D}$ be given by \eqref{charafunchi}. Then, for any integer $n$ and nonnegative real numbers $r,s$, $m=\max\{r,s\}$, it holds
\begin{multline}
 \frac{a}{c}\mathcal{D}\bigg(\frac{ce}{a},\frac{bc}{a},\frac{cd}{a},\frac{b c d e}{a}q^{2-4r}\bigg)\sum_{k=0}^{n-1}\frac{1-
 \frac{b c d e}{a}q^{4rk+2-4r}}{1-
 \frac{b c d e}{a}q^{2-4r}}
\frac{(bq,dq,eq,bc^2deq/a^2;q^{2r})_{k-1}}
{(cdeq/a,bceq/a,bcdq/a,aq/c;q^{2r})_{k}}\\
\times\frac{(yq,wq,uq,yz^2wuq/x^2;q^{2s})_{k}}
{(zwuq/x,yzuq/x,yzwq/x,xq/z;q^{2s})_{k}} q^{2rk+1-2r}\\
=
\frac{(bq,dq,eq,bc^2deq/a^2;q^{2r})_{n-1}}
{(cdeq/a,bceq/a,bcdq/a,aq/c;q^{2r})_{n-1}}
\frac{(yq,wq,uq,yz^2wuq/x^2;q^{2s})_{n}}
{(zwuq/x,yzuq/x,yzwq/x,xq/z;q^{2s})_{n}}\\
-\frac{(q^{2r-1}-cde/a)(q^{2r-1}-bce/a)(q^{2r-1}-bcd/a)(q^{2r-1}-a/c)}{(q^{2r-1}-b)(q^{2r-1}-d)(q^{2r-1}-e)(q^{2r-1}-bc^2de/a^2)}\\
-\frac{x}{z}\mathcal{D}\bigg(\frac{yz}{x},\frac{zu}{x},\frac{zw}{x},\frac{y z w u}{x} q^{2}\bigg)\sum_{k=0}^{n-1}\frac{1-\frac{y z w u}{x} q^{4s k+2}}{1-\frac{y z w u}{x} q^2}
\frac{(yq,wq,uq,yz^2wuq/x^2;q^{2s})_{k}}
{(zwuq/x,yzuq/x,yzwq/x,xq/z;q^{2s})_{k+1}}\label{pppppp-generalid-new}\\\times
\frac{(bq,dq,eq,bc^2deq/a^2;q^{2r})_{k}}
{(cdeq/a,bceq/a,bcdq/a,aq/c;q^{2r})_{k}} q^{2sk+1}.
 \end{multline}
 \end{tl}
\pf It suffices to set in \eqref{pppppp-generalid} that $p=q, r_i=r,  s_i=s,~1\leq i\leq 4$.\qed

\begin{tl}[$(2,2)$-type transformation with degree $2m$] \label{degree-2s2-1}
For any nonzero complex numbers $b,d,e,y,u,w,q$ and $m=\max\{r_i,s_i|1\leq i\leq 2\}$ and $t_1,t_2>0$, define
\begin{align}
C_0(t_1)=(1-bt_1)(1-dt_1)/t_1,~~D_0(t_2)=(1-yt_2)(1-wt_2)/t_2,
\label{xxxxx-new}
\end{align}
 there holds
\begin{multline}
C_0(t_1)\sum_{k=0}^{n-1}(1-et_1q^{2(r_2-r_1)(k-1)})
\big(1-b d et_1q^{2(r_1+r_2)(k-1)+2} \big)q^{2r_1(k-1)+1}\\
\times
\frac{(bq,dq;q^{2r_1})_{k-1}(eq,bdet_1^2q;q^{2r_2})_{k-1}}
{(det_1q,bet_1q;q^{2r_2})_{k}(bdt_1q, q/t_1;q^{2r_1})_{k}}\times\frac{(yq,wq;q^{2s_1})_{k}(uq,ywut_2^2q;q^{2s_2})_{k}}
{(wut_2q,yut_2q;q^{2s_2})_{k}(ywt_2q, q/t_2;q^{2s_1})_{k}}\\
=
\frac{(bq,dq;q^{2r_1})_{n-1}(eq, bdet_1^2q;q^{2r_2})_{n-1}}
{(det_1q,bet_1q;q^{2r_2})_{n-1}(bdt_1q, q/t_1;q^{2r_1})_{n-1}}\times\frac{(yq,wq;q^{2s_1})_{n}(uq,ywut_2^2q;q^{2s_2})_{n}}
{(wut_2q,yut_2q;q^{2s_2})_{n}(ywt_2q, q/t_2;q^{2s_1})_{n}}\\
-\frac
{(q^{2r_2}-det_1)(q^{2r_2}-bet_1)(q^{2r_1}-bdt_1)(q^{2r_1}-q/t_1)}{(q^{2r_1}-bq)(q^{2r_1}-dq)(q^{2r_2}-eq)(q^{2r_2}-bdet_1^2q)}\\
-D_0(t_2)\sum_{k=0}^{n-1}(1-ut_2q^{2(s_2-s_1)k})\big(1-yw ut_2 q^{2(s_1+s_2)k+2}\big)q^{2s_1k+1}\label{4.4}\\
\times\frac{(yq,wq;q^{2s_1})_{k}(uq,ywut_2^2q;q^{2s_2})_{k}}
{(wut_2q,yut_2q;q^{2s_2})_{k+1}(ywt_2q,q/t_2;q^{2s_1})_{k+1}}\times
\frac{(bq,dq;q^{2r_1})_{k})(eq,bdet_1^2q;q^{2r_2})_{k}}
{(det_1q,bet_1q;q^{2r_2})_{k}(bdt_1q,q/t_1;q^{2r_1})_{k}}.
\end{multline}
\end{tl}
\pf It is sufficient to consider $p=q$ and then make the substitutions   $
(r_1,r_2,r_3,r_4)\to(r_1,r_1,r_2,r_2)$, $(s_1,s_2,s_3,s_4)\to(s_1,s_1,s_2,s_2)$
in Theorem \ref{theorem41}. Finally,  making the substitutions   $c=at_1, z=xt_2$ and use
$C_0(t_1)$ and $D_0(t_2)$ given by \eqref{xxxxx-new}, we obtain \eqref{4.4} .
\qed

\begin{remark}In this view, we point out that  the quadratic transformation \eqref{ooooo} is  also a $(2,2)$-type with degree $2$. To make this clear, we need only to take substitutions in
Corollary \ref{degree-2s2-1}
$
(b,d,e,t_1)\to (d,q/d,bcq^2,a/q)$, $(y,w,u,t_2)\to (abc/q,q/a,ac,1/c)$, and
 $(r_1,r_2,s_1,s_2)=(1/2,1,1/2,1)$.
\end{remark}

\begin{tl}[$(2,2)$-type transformation with degree $2m$] \label{degree-2s2-3}
For any nonzero complex numbers $b,d,e,y,u,w,q$ and $t_1,t_2$, $C_0(t_1)$ be the same as \eqref{xxxxx-new},
 there holds
\begin{multline}
C_0(t_1)\sum_{k=0}^{n-1}\mathcal{D}\big(et_1q^{2(r_2-r_1)(k-1)},b d et_1q^{2(r_1+r_2)(k-1)+2} \big)q^{2r_1(k-1)+1}\\
\times
\frac{(bq,dq;q^{2r_1})_{k-1}(eq,bdet_1^2q;q^{2r_2})_{k-1}}
{(det_1q,bet_1q;q^{2r_2})_{k}(bdt_1q, q/t_1;q^{2r_1})_{k}}\times\frac{(yq;q^{2s_1})_{k}(wq,uq,ywut_2^2q;q^{2s_2})_{k}}
{(wut_2q;q^{3s_2-s_1})_{k}(yut_2q,ywt_2q, q/t_2;q^{s_1+s_2})_{k}}\\
=
\frac{(bq,dq;q^{2r_1})_{n-1}(eq, bdet_1^2q;q^{2r_2})_{n-1}}
{(det_1q,bet_1q;q^{2r_2})_{n-1}(bdt_1q, q/t_1;q^{2r_1})_{n-1}}\times\frac{(yq;q^{2s_1})_{n}(wq,uq,ywut_2^2q;q^{2s_2})_{n}}
{(wut_2q;q^{3s_2-s_1})_{n}(yut_2q,ywt_2q, q/t_2;q^{s_1+s_2})_{n}}\\
-\frac
{(q^{2r_2}-det_1q)(q^{2r_2}-bet_1q)(q^{2r_1}-bdt_1q)(q^{2r_1}-q/t_1)}{(q^{2r_1}-bq)(q^{2r_1}-dq)(q^{2r_2}-eq)(q^{2r_2}-bdet_1^2q)}\\
-\frac{1}{t_2}\sum_{k=0}^{n-1}\mathcal{D}\bigg(yt_2q^{(s_1-s_2)k},ut_2q^{(s_2-s_1)k},wt_2q^{(s_2-s_1)k},yw ut_2 q^{(s_1+3s_2)k+2}\bigg)q^{(s_1+s_2)k+1}\\
\times\frac{(yq;q^{2s_1})_{k}(wq,uq,ywut_2^2q;q^{2s_2})_{k}}
{(wut_2q;q^{3s_2-s_1})_{k+1}(yut_2q,ywt_2q,q/t_2;q^{s_1+s_2})_{k+1}}\times
\frac{(bq,dq;q^{2r_1})_{k})(eq,bdet_1^2q;q^{2r_2})_{k}}
{(det_1q,bet_1q;q^{2r_2})_{k}(bdt_1q,q/t_1;q^{2r_1})_{k}}.
\end{multline}
\end{tl}
\pf It is sufficient to make the substitutions $p=q;
(r_1,r_2,r_3,r_4)\to(r_1,r_1,r_2,r_2)$,~$(s_1,s_2,s_3,s_4)\\\to(s_1,s_2,s_2,s_2)$
in Theorem \ref{theorem41} and then replace  $(c,z)$ with $(at_1,xt_2)$ and use
$C_0(t_1)$ given by \eqref{xxxxx-new}.
\qed

\begin{remark}
As a matter of fact,
 the quartic transformation  \eqref{wantedresult} is just  a $(2,2)$-type with degree $4$,
 if we make the substitutions in Corollary \ref{degree-2s2-3} $
(b,d,e,t_1)\to (a,b,a^2b^2q,1/ab)$ and $(y,w,u,t_2)\to (ab,ab/q,ab/q^2,1/ab)$, as well as let $(r_1,r_2,s_1,s_2)=(1/2,2,1/2,3/2)$.
\end{remark}

\begin{tl}[$(2,3)$-type transformation with degree $2m$]  \label{degree-2s2-2} Let $\mathcal{D}$ be given by \eqref{charafunchi}.  Then, for any nonzero complex numbers $b,d,e,y,u,w,q, t_1,t_2$, and  $m=\max\{r_i,s_j|i=1,2,j=1,2,3\}$, 
there holds
\begin{multline}
C_0(t_1)\sum_{k=0}^{n-1}\mathcal{D}\big(et_1q^{2(r_2-r_1)(k-1)},bdet_1q^{2(r_1+r_2)(k-1)+2}\big) q^{2r_1(k-1)+1}
\times
\frac{(bq,dq;q^{2r_1})_{k-1}(eq,bdet_1^2q;q^{2r_2})_{k-1}}
{(det_1q,bet_1q;q^{2r_2})_{k}(bdt_1q,q/t_1;q^{2r_1})_{k}}\\
\times
\frac{(yq;q^{2s_1})_{k}(wq,uq;q^{2s_2})_{k}(ywut_2^2q;q^{2s_3})_k}
{(wut_2q;q^{2s_2+s_3-s_1})_k(yut_2q,ywt_2q;q^{s_1+s_3})_k (q/t_2;q^{s_1+2s_2-s_3})_{k}}\\
=
\frac{(bq,dq;q^{2r_1})_{n-1}(eq,bdet_1^2q;q^{2r_2})_{n-1}}
{(det_1q,bet_1q;q^{2r_2})_{n-1}(bdt_1q,q/t_1;q^{2r_1})_{n-1}}
\times
\frac{(yq;q^{2s_1})_{n}(wq,uq;q^{2s_2})_{n}(ywut_2^2q;q^{2s_3})_n}
{(wut_2q;q^{2s_2+s_3-s_1})_n(yut_2q,ywt_2q;q^{s_1+s_3})_n (q/t_2;q^{s_1+2s_2-s_3})_{n}}\\
-\frac
{(q^{2r_2}-det_1q)(q^{2r_2}-bet_1q)(q^{2r_1}-bdt_1q)(q^{2r_1}-q/t_1)}
{(q^{2r_1}-bq)(q^{2r_1}-dq)(q^{2r_2}-eq)(q^{2r_2}-bdet_1^2q)}
\end{multline}
 \begin{multline}
-\frac{1}{t_2}\sum_{k=0}^{n-1}\mathcal{D}\bigg(yt_2q^{(s_1+s_3-2s_2)k},ut_2q^{(s_3-s_1)k},wt_2q^{(s_3-s_1)k},ywut_2q^{(s_1+2s_2+s_3)k+2}\bigg)q^{(s_1+2s_2-s_3)k+1}
\nonumber\\
\times
\frac{(yq;q^{2s_1})_{k}(wq,uq;q^{2s_2})_{k}(ywut_2^2q;q^{2s_3})_k}
{(wut_2q;q^{2s_2+s_3-s_1})_{k+1}(yut_2q,ywt_2q;q^{s_1+s_3})_{k+1} (q/t_2;q^{s_1+2s_2-s_3})_{k+1}}
\times
\frac{(bq,dq;q^{2r_1})_{k}(eq,bdet_1^2q;q^{2r_2})_{k}}
{(det_1q,bet_1q;q^{2r_2})_{k}(bdt_1q,q/t_1;q^{2r_1})_{k}}.
\end{multline}
\end{tl}
\pf It is sufficient to make the substitutions $p=q$,  $
(r_1,r_2,r_3,r_4)\to(r_1,r_1,r_2,r_2)$ and $(s_1,s_2,s_3,s_4)\to(s_1,s_2,s_2,s_3)$
in Theorem \ref{theorem41} and then make the parametric replacement  $c=at_1, z=xt_2$
and use
$C_0(t_1)$
given by \eqref{xxxxx-new}
.
\qed

\begin{remark}Indeed, the cubic transformation \eqref{pppppp-revised}  turns out to be $(2,2)$-type with degree $3$. To see this, we need to
make the substitutions in Corollary \ref{degree-2s2-2}
$
(b,d,e,t_1)\to (a,1/ab,q^2c,ab/c)$ and $
(y,w,u,t_2)\to (b/q,1/b,q/b,qab)
$,  $(r_1,r_2)=(1/2,3/2)$, $(s_1,s_2,s_3)=(1/2,1,3/2)$.
\end{remark}

\begin{tl}[$(3,3)$-type transformation with degree $2m$]  \label{degree-3}
For any nonzero complex numbers $b,d,e,y,u,w,q$ and $r_i,s_i,i=1,2,3$,  there holds
 \begin{multline}
\frac{1}{t_1}\sum_{k=0}^{n-1}\mathcal{D}\bigg(et_1q^{(r_2+r_3-2r_1)(k-1)},bt_1q^{(r_3-r_2)(k-1)},
 dt_1q^{(r_3-r_2)(k-1)},bdet_1q^{r(k-1)+2} \big)\\
\times q^{(r-2r_3)(k-1)+1}
\frac{(bq,dq;q^{2r_1})_{k-1}(eq;q^{2r_2})_{k-1}
(bdet_1^2q;q^{2r_3})_{k-1}}
{(deqt_1,beqt_1;q^{r-2r_1})_{k}(bdqt_1;q^{r-2r_2})_{k}
(q/t_1;q^{r-2r_3})_{k}}\\
=
\frac{(bq,dq;q^{2r_1})_{n-1}(eq;q^{2r_2})_{n-1}
(bt_1^2deq;q^{2r_3})_{n-1}}
{(deqt_1,beqt_1;q^{r-2r_1})_{n-1}(bdqt_1;q^{r-2r_2})_{n-1}
(q/t_1;q^{r-2r_3})_{n-1}}
\times\frac{(yq,wq;q^{2s_1})_{n}(uq;q^{2s_2})_{n}
(yt_2^2wuq;q^{2s_3})_{n}}
{(wuqt_2,yuqt_2;q^{s-2s_1})_{n}
(ywqt_2;q^{s-2s_2})_{n}(q/t_2;q^{s-2s_3})_{n}}
\\
-\frac
{(q^{r-2r_1}-deqt_1)(q^{r-2r_1}-beqt_1)(q^{r-2r_2}-bdqt_1)(q^{r-2r_3}-q/t_1)}{(q^{2r_1}-bq)(q^{2r_1}-dq)(q^{2r_2}-eq)(q^{2r_3}-bt_1^2deq)}\\
-\frac{1}{t_2}\sum_{k=0}^{n-1}\mathcal{D}\bigg(yt_2q^{(s_3-s_2)k},ut_2q^{(s_2+s_3-2s_1)k},wt_2q^{(s_3-s_2)k},ywut_2q^{sk+2}\bigg)q^{(s-2s_3)k+1}\\
\times\frac{(yq,wq;q^{2s_1})_{k}(uq;q^{2s_2})_{k}
(yt_2^2wuq;q^{2s_3})_{k}}
{(wuqt_2,yuqt_2;q^{s-2s_1})_{k+1}
(ywqt_2;q^{s-2s_2})_{k+1}(q/t_2;q^{s-2s_3})_{k}}\\
\times
\frac{(bq,dq;q^{2r_1})_{k}(eq;q^{2r_2})_{k}
(bdeqt_1^2;q^{2r_3})_{k}}
{(deqt_1,beqt_1;q^{r-2r_1})_{k}(bdqt_1;q^{r-2r_2})_{k}
(q/t_1;q^{r-2r_3})_{k}}.
 \end{multline}
\end{tl}
\pf It is sufficient to make the substitutions $p=q$,$
(r_1,r_2,r_3,r_4)\to(r_1,r_1,r_2,r_3)$, $(s_1,s_2,s_3,s_4)\to(s_1,s_1,s_2,s_3)$ and then $c=at_1, z=xt_2$
in Theorem \ref{theorem41}.
\qed

We end our paper with Chu's recent result which can be regarded as a $(3,3)$-type transformation with degree $5$, i.e., a special case of  Corollary \ref{degree-3}.

\begin{lz}[{\rm\!\! Quintic transformation: \cite[Sec.2.1]{chulatest}}]Let $\mathcal{D}$ be the same as in Lemma \ref{lemm12}, and for any integer $n$, we define
\begin{align}
G_n(a, c):=\sum_{k=0}^n\left(1- a^3q^{6 k-1}\right)
\frac{\left(a^2 c, a^3 / c ; q^5\right)_k}{(c, a / c ; q)_k}
\frac{(q a, a / q ; q)_{2 k}}{\left(q a^2 ; q\right)_{4 k}} q^k.\label{chufun}
\end{align}
Then
\begin{align}
G_n(a, c)&=G_n\left(q a, q^3 c\right)\frac{\mathcal{D}(q a, a, a / q, a^2 c)}{\mathcal{D}(a^{-1}, c, q c, q^2 c)} \frac{(q c / a ; q)_2}{\left(q a^2 ; q\right)_2}+\frac{\mathcal{D}(a c, a^2, a^2 / q, q c / a)}{\mathcal{D}(a,c,cq)} \label{4.6}\\
&+ \frac{\mathcal{D}(q^{3+3 n} a c, q a, 1/a, a / q, a^2 c, a^3 / c )}{\mathcal{D}(c, 1/q c, q^2 c, q a^2, q^2 a^2)} \frac{\left(q^5 a^2 c, q^5 a^3 / c ; q^5\right)_n}{\left(q^3 c, a / c ; q\right)_n}
\frac{\left(q a, q^2 a ; q\right)_{2 n}}{\left(q^3 a^2 ; q\right)_{4 n}} .\nonumber
\end{align}
\end{lz}
\pf To show \eqref{4.6},we need to make the substitutions in
{Corollary \ref{degree-3} that $(b,d,e,t_1)\to (a,aq,a^2q/c,cq/a)$ and $
(y,w,u,t_2)\to (a/q ,a,acq^2,1/cq).$
Then we have
 \begin{multline}
\frac{a}{c}\sum_{k=0}^{n-1}\mathcal{D}\big(aq^{(r_2+r_3-2r_1)(k-1)+2},cq^{(r_3-r_2)(k-1)+1},cq^{(r_3-r_2)(k-1)+2},a^3q^{r(k-1)+5} \big)q^{(r-2r_3)(k-1)}\nonumber\\
\times\frac{(aq,aq^2;q^{2r_1})_{k-1}(a^2q^2/c;q^{2r_2})_{k-1}
(a^2cq^5;q^{2r_3})_{k-1}}
{(a^2 q^4,a^2 q^3;q^{r-2r_1})_{k}(a c q^3;q^{r-2r_2})_{k}
(a/c;q^{r-2r_3})_{k}}\times\frac{(a,aq;q^{2s_1})_{k}(acq^3;q^{2s_2})_{k}
(a^3/c;q^{2s_3})_{k}}
{(a^2 q^2,a^2 q;q^{s-2s_1})_{k}
(a^2/c q;q^{s-2s_2})_{k}(c q^2;q^{s-2s_3})_{k}}\nonumber\\
=\frac{(a q,a q^2;q^{2r_1})_{n-1}(a^2q^2/c;q^{2r_2})_{n-1}
(a^2 c q^5;q^{2r_3})_{n-1}}{(a^2 q^4,a^2 q^3;q^{r-2r_1})_{n-1}(a c q^3;q^{r-2r_2})_{n-1}
(a/c;q^{r-2r_3})_{n-1}}\times\frac{(a,aq;q^{2s_1})_{n}(acq^3;q^{2s_2})_{n}
(a^3/c;q^{2s_3})_{n}}
{(a^2 q^2,a^2 q;q^{s-2s_1})_{n}
(a^2/c q;q^{s-2s_2})_{n}(cq^2;q^{s-2s_3})_{n}}
\\
-\frac{(q^{r-2r_1}-a^2 q^4)(q^{r-2r_1}-a^2 q^3)(q^{r-2r_2}-a c q^3)(q^{r-2r_3}-a/c)}{(q^{2r_1}-a q)(q^{2r_1}-a q^2)(q^{2r_2}-a^2q^2/c)(q^{2r_3}-a^2 c q^5)}\nonumber\\
-cq^2\sum_{k=0}^{n-1}\mathcal{D}\big(aq^{(s_3-s_2)k-2}/c,aq^{(s_2+s_3-2s_1)k+1},aq^{(s_3-s_2)k-1}/c,a^3q^{s k+2}\big)q^{(s-2s_3)k}\\
\times\frac{(a,aq;q^{2s_1})_{k}(acq^3;q^{2s_2})_{k}
(a^3/c;q^{2s_3})_{k}}
{(a^2 q^2,a^2 q;q^{s-2s_1})_{k+1}
(a^2/c q;q^{s-2s_2})_{k+1}(c q^2;q^{s-2s_3})_{k+1}}\times\frac{(a q,a q^2;q^{2r_1})_{k}(a^2q^2/c;q^{2r_2})_{k}
(a^2 c q^5;q^{2r_3})_{k}}
{(a^2 q^4,a^2 q^3;q^{r-2r_1})_{k}(a c q^3;q^{r-2r_2})_{k}
(a/c;q^{r-2r_3})_{k}}.
 \end{multline}
From here, on taking $(r_1,r_2,r_3)=(s_1,s_2,s_3)=(1,3/2,5/2)$, we come up with  \begin{multline}
\frac{a}{c}\sum_{k=0}^{n-1}\mathcal{D}\big(aq^{2k},cq^{k},cq^{k+1},a^3q^{6k-1}\big)q^{k-1}\nonumber\\
\times\frac{(aq,aq^2;q^{2})_{k-1}(a^2q^2/c;q^{3})_{k-1}
(a^2cq^5;q^{5})_{k-1}}
{(a^2 q^4,a^2 q^3;q^{r-2})_{k}(a c q^3;q^{3})_{k}
(a/c;q)_{k}}
\times\frac{(a,aq;q^{2})_{k}(acq^3;q^{3})_{k}
(a^3/c;q^{5})_{k}}
{(a^2 q^2,a^2 q;q^{4})_{k}
(a^2/c q;q^{3})_{k}(c q^2;q)_{k}}\nonumber\\
=\frac{(a q,a q^2;q^{2})_{n-1}(a^2 q^2/c;q^{3})_{n-1}
(a^2 c q^5;q^{5})_{n-1}}{(a^2 q^4,a^2 q^3;q^{4})_{n-1}(a c q^3;q^{3})_{n-1}
(a/c;q)_{n-1}}
\times\frac{(a,aq;q^{2})_{n}(acq^3;q^{3})_{n}
(a^3/c;q^{5})_{n}}
{(a^2 q^2,a^2 q;q^{4})_{n}
(a^2/c q;q^{3})_{n}(cq^2;q)_{n}}\\
-\frac{(q^{4}-a^2 q^4)(q^{4}-a^2 q^3)(q^{3}-a c q^3)(q-a/c)}{(q^{2}-a q)(q^{2}-a q^2)(q^{3}-a^2 q^2/c)(q^{5}-a^2 c q^5)}\\-cq^2\sum_{k=0}^{n-1}\mathcal{D}\big(aq^{k-2}/c,aq^{2k+1},aq^{k-1}/c,a^3q^{6k+2}\big)q^{k}\\
\times\frac{(a,aq;q^{2})_{k}(acq^3;q^{3})_{k}
(a^3/c;q^{5})_{k}}
{(a^2 q^2,a^2 q;q^{4})_{k+1}
(a^2/c q;q^{3})_{k+1}(c q^2;q)_{k+1}}
\times\frac{(a q,a q^2;q^{2})_{k}(a^2 q^2/c;q^{3})_{k}
(a^2 c q^5;q^{5})_{k}}
{(a^2 q^4,a^2 q^3;q^{4})_{k}(a c q^3;q^{3})_{k}
(a/c;q)_{k}}.\nonumber
 \end{multline}
That is
 \begin{multline}
C_0+\frac{a}{c}\sum_{k=2}^{n-1}\big(1-aq^{2k}\big)
\left(1-a^3q^{6k-1} \right)q^{k-1}\nonumber\\
\times\frac{(aq,aq^2;q^{2})_{k-1}(a^2q^2/c;q^{3})_{k-1}
(a^2cq^5;q^{5})_{k-1}}
{(a^2 q^4,a^2 q^3;q^{4})_{k}(a c q^3;q^{3})_{k}
(a/c;q)_{k}}
\times\frac{(a,aq;q^{2})_{k}(acq^3;q^{3})_{k}
(a^3/c;q^{5})_{k}}
{(a^2 q^2,a^2 q;q^{4})_{k}
(a^2/c q;q^{3})_{k}(c q^2;q)_{k-2}}\nonumber\\
=\frac{(a q,a q^2;q^{2})_{n-1}(a^2 q^2/c;q^{3})_{n-1}
(a^2 c q^5;q^{5})_{n-1}}{(a^2 q^4,a^2 q^3;q^{4})_{n-1}(a c q^3;q^{3})_{n-1}
(a/c;q)_{n-1}}\times\frac{(a,aq;q^{2})_{n}(acq^3;q^{3})_{n}
(a^3/c;q^{5})_{n}}
{(a^2 q^2,a^2 q;q^{4})_{n}
(a^2/c q;q^{3})_{n}(cq^2;q)_{n}}\\
-\frac{(q^{4}-a^2 q^4)(q^{4}-a^2 q^3)(q^{3}-a c q^3)(q-a/c)}{(q^{2}-a q)(q^{2}-a q^2)(q^{3}-a^2q^2/c)(q^{5}-a^2 c q^5)}-D_0
-cq^2\sum_{k=2}^{n-1}\big(1-aq^{2k+1}\big)\left(1-a^3q^{6k+2}\right)q^{k}\\
\times\frac{(a,aq;q^{2})_{k}(acq^3;q^{3})_{k}
(a^3/c;q^{5})_{k}}
{(a^2 q^2,a^2 q;q^{4})_{k+1}
(a^2/c q;q^{3})_{k+1}(c q^2;q)_{k+1}}\times\frac{(a q,a q^2;q^{2})_{k}(a^2 q^2/c;q^{3})_{k}
(a^2 c q^5;q^{5})_{k}}
{(a^2 q^4,a^2 q^3;q^{4})_{k}(a c q^3;q^{3})_{k}
(a/c;q)_{k-2}}.
 \end{multline}
Referring  to \eqref{chufun}, we obtain the conclusion as desired.
\qed

\section*{Acknowledgements}The authors thank   the anonymous referees  for their  careful reading and  insightful comments on the first version of this manuscript.
This  work was supported by the National Natural Science Foundation of
China [Grant No. 11971341].


\end{document}